\documentclass{amsart}
\usepackage{graphicx} 
\usepackage{preamble}
\usepackage{spectralsequences} 
\usepackage{subcaption} 
 
\title{Spoke topological Hochschild homology}
\author{Gabriel Angelini-Knoll}
\address{Department of Mathematics, Applied Mathematics, and Statistics, Case
Western Reserve University, Cleveland, OH, USA}
\email{gabriel.angelini-knoll@case.edu}

\author{Foling Zou}
\address{State Key Laboratory of Mathematical Sciences, Academy of Mathematics and Systems Science, Chinese Academy of Sciences, Beijing 100190, China}
\email{zoufoling@amss.ac.cn}
\date{}

\begin{document}
\begin{abstract} 
Fix primes $p$ and $\ell$, and let $C_p$ be the cyclic group of order $p$. We compute the $C_p$-equivariant spoke topological Hochschild homology of $\mF_{\ell}$ and prove it exhibits a form of B\"okstedt periodicity. Here spoke topological Hochschild homology is a variant of topological Hochschild homology where one replaces the circle in the construction with the unreduced suspension of $C_p$.
As an application, we use this result to give a new proof of the Segal conjecture for the cyclic group of order an odd prime $p$. 
\end{abstract}

\maketitle 

\sseqset{ axes type = center, axes gap =0cm, xscale = 0.75, yscale = 1.6}


\DeclareSseqGroup\towerGG {} {
\class(0,0)
\foreach \y in {1,2} {
\class(0,-\y)
\structline
}
\class(2,-1)
\structline(0,0,-1)
\class(4,-2)
\structline

\class(1,-2) \class(3,-2)
\class(1,-1)
\structline(1,-2,-1) \structline(3,-2,-1)
\class(2,-2)
\structline(0,-1,-1) \structline(2,-1,-1)
}

\DeclareSseqGroup\towerGGTwice {} {
\class(0,0)
\foreach \y in {2,4} {
\class(0,-\y)
\structline
}
\class(2,-2)
\structline(0,0,-1)
\class(4,-4)
\structline

\class(1,-4) \class(3,-4)
\class(1,-2)
\structline(1,-4,-1) \structline(3,-4,-1)
\class(2,-4)
\structline(0,-2,-1) \structline(2,-2,-1)
}

\DeclareSseqGroup\trap{}{
  \class(0,0)
  \class[blue](0,-1)
  \structline
  \class(1,-2)
  \structline[blue]
  \class[blue](1,-1)
  \structline
  \structline[blue](0,0,-1)
}
\DeclareSseqGroup\towerGSpokeGrd{}{
  \trap(0,0)
  \trap(0,-2)
  \trap(2,-2)
  \foreach \x in {0,1,2}{
    \class(2*\x,-4)}
  \structline(0,-4,-1)(0,-3,-1)
  \foreach \x in {0,1,2}{
    \structline(\x,-2-\x,-1)(\x,-1-\x,-1)}
}

\DeclareSseqGroup\towerBB {} {
\class(0,0)
\foreach \y in {1,2} {
\class(0,\y)
\structline
}
\class(-2,1)
\structline(0,0,-1)
\class(-4,2)
\structline

\class(-1,2) \class(-3,2)
\class(-1,1)
\structline(-1,2,-1) \structline(-3,2,-1)
\class(-2,2)
\structline(0,1,-1) \structline(-2,1,-1)
}
\DeclareSseqGroup\towerG {} {
\class(0,0)
\class(0,-1)
\structline
\class(2,-1)
\structline(0,0,-1)
\class(1,-1)
}
\DeclareSseqGroup\towerB {} {
\class(0,0)
\class(0,1)
\structline
\class(-2,1)
\structline(0,0,-1)
\class(-1,1)
}
\DeclareSseqGroup\towerGtiny {} {
\class(0,0)
\class(0,-1)
\structline
\class(2,-1)
\structline(0,0,-1)
}
\DeclareSseqGroup\towerBtiny {} {
\class(0,0)
\class(0,1)
\structline
\class(-2,1)
\structline(0,0,-1)
}

\DeclareSseqGroup\towerSpoke {} {
  \foreach \x in {0,1}{
    \class(\x,-1)
    \class(\x+2,-1)
    \class(\x,0)
    \structline(\x,-1,-1)
    \structline(\x+2,-1,-1)
    }
  }

  \DeclareSseqGroup\towerSpokeTwice{}{
  \foreach \x in {0,1}{
    \class(\x,-2)
    \class(\x+2,-2)
    \class(\x,0)
    \structline(\x,-2,-1)
    \structline(\x+2,-2,-1)
    }}
\DeclareSseqGroup\towerMM {} {
\class(1,-1)
\class(1,-2)
\structline
\class(2,-1)
\structline(1,-1,-1)
\class(4,-2)
\structline

\class(-1,1)
\class(-1,2)
\structline
\class(-2,1)
\structline(-1,1,-1)
\class(-4,2)
\structline
}
\DeclareSseqGroup\towerPhi {} {
  \class(0,0)
  \DoUntilOutOfBounds{
    \class(\lastx,\lasty+1)
    \structline
  }
  \class(0,-1) \structline(0,0,-1)
  \DoUntilOutOfBounds{
    \class(\lastx,\lasty-1)
    \structline
  }
}

\begin{sseqdata}[ name = Hstar,
  title = {},
  classes = { circle, fill },
  class labels = { above left = 0.3em , pin},
  x range = {0}{4},
        y tick handler = {
    \ifnum#1 = 0\relax
      \vphantom{0}
    \else
      \ifnum#1 = 1 \relax
        \protect\vphantom{2}\lambda
      \else
        \ifnum#1 = -1 \relax
          \protect\vphantom{2}-\lambda
        \else
          #1\lambda
        \fi
      \fi
    \fi}
  ]

\towerGG[tag = x](0,0)

\classoptions["a_{\lambda}"{above left}](0,-1)
\classoptions["{\kappa_{\lambda}}"{above left}](1,-1)
\classoptions["{u_{\lambda}}"{above}](2,-1)

\end{sseqdata}

\begin{sseqdata}[ name = HstarSpoke, title = {},
  classes = { circle, fill }, y range = {-2}{0},
  y tick handler = {
    \ifnum#1 = 0\relax
      \vphantom{0}
    \else
      \ifnum#1 = 1 \relax
        \protect\vphantom{2}\lambda
      \else
        \ifnum#1 = -1 \relax
          \protect\vphantom{2}-\lambda
        \else
          #1\lambda
        \fi
      \fi
    \fi}
  ]
\towerSpoke[blue](0,-1)
\end{sseqdata}

\begin{sseqdata}[ name = HstarFull, title = {},
  classes = { circle, fill }, yscale = 0.5,
  y tick handler = {
    \ifnum#1 = 0\relax
      \vphantom{0}
    \else
      \ifnum#1 = 1 \relax
        \protect\vphantom{2}\Yright
      \else
        \ifnum#1 = -1 \relax
          \protect\vphantom{2}-\Yright
        \else
           #1\Yright
        \fi
      \fi
    \fi}
  ]
  \towerGSpokeGrd(0,0)
\end{sseqdata}

\begin{sseqdata}[ name = HstarFullShift, title = {},
  classes = { circle, fill }, yscale = 0.5,
  y range = {-5}{0},
  y tick handler = {
    \ifnum#1 = 0\relax
      \vphantom{0}
    \else
      \ifnum#1 = 1 \relax
        \protect\vphantom{2}\Yright
      \else
        \ifnum#1 = -1 \relax
          \protect\vphantom{2}-\Yright
        \else
           #1\Yright
        \fi
      \fi
    \fi}
  ]
  \towerGSpokeGrd(0,-1)

\class(0,-1)
\foreach \x in {1,2}{
    \class(2*\x,-2*\x-1)
    \structline
    }
\end{sseqdata}

\begin{sseqdata}[ name = HTstar, title = {},
   classes = { circle, fill },
   class labels = { above left = 0.3em , pin},
   y range = {0}{2},
   x range = {0}{5},
  y tick handler = {
    \ifnum#1 = 0\relax
      \vphantom{0}
    \else
      \ifnum#1 = 1 \relax
        \protect\vphantom{2}\lambda
      \else
        \ifnum#1 = -1 \relax
          \protect\vphantom{2}-\lambda
        \else
          #1\lambda
        \fi
      \fi
    \fi}
  ]

  \towerSpoke(1,1)
  \classoptions["{\widehat{\xi}_1}"{above left}](1,1)
  \classoptions["{\underline{\xi}_1}"{above left}](2,1)
  
  \towerSpoke(2,1)
  \classoptions["{\widehat{\tau}_1}"{above}](2,1,-1)
  \classoptions["{\underline{\tau}_1}"{above}](3,1)
\end{sseqdata}








\section{Introduction}
Topological Hochschild homology (THH) has a geometric interpretation as the tensoring with the circle
in commutative ring spectra, by an insightful theorem of
McClure--Schw\"anzl--Vogt~\cite{MSV97}. It has applications to string topology, deformation theory of $\mathbb{E}_1$ algebras, and algebraic K-theory. 
The fact that the topological Hochschild homology of a finite field $\mathbb{F}_p$ of order $p$ is a polynomial algebra on a generator in degree two, a result sometimes called B\"okstedt periodicity, is a computation of fundamental importance~\cite{Bok85}. For example, B\"okstedt periodicity can be used to prove Bott periodicity~\cite{HN20}. 

In the $C_p$-equivariant setting, one may replace the circle with the spoke circle $S^{\spoke}$ defined as the unreduced suspension of $C_p$. Here $C_p$ denotes a cyclic group of order $p$ where $p$ is a prime. Using the equivariant Loday construction of~\cite{LRZLoday}, one can make sense of tensoring a $C_p$-commutative ring spectrum $R$ with the spoke sphere, a construction we call spoke topological Hochschild homology and denote $\THH^{\spoke}(R)$. In particular, spoke topological Hochschild homology $\THH^{\spoke}(R)$ can be identified with $R\otimes_{N_e^{C_p}i_e^*R}R$. When $p=2$, this produces Real topological Hochschild homology~\cite{DMPR21,QS21}, which has applications to the study of quadratic forms and surgery theory of manifolds. When $p$ is an odd prime, it is expected to have applications to a $C_p$-equivariant analogue of free loop space theory, deformations of $\mathbb{E}_{\spoke}$ algebras, algebraic K-theory, and to computations of $C_p$-fixed points of $C_p$-equivariant norms.

Let $p$ and $\ell$ be primes. Let $\mF_\ell$ denote the Eilenberg--MacLane spectrum of the constant $C_p$-equivariant Mackey functor associated to a finite field of order $\ell$. Our first main result is the analogue of B\"okstedt periodicity for this new theory.  
\begin{thm}[B\"okstedt periodicity for spoke THH]\label{thm:Bok}
There is an equivalence
\[ 
\mF_\ell[S^{1+\spoke}]\simeq \sTHH(\mF_\ell) \,.
\]
\end{thm}
\noindent Here, for a $C_p$-commutative ring spectum $R$,
$R[S^{1+\spoke}]:=R\otimes \left ( N_e^{C_p}\bS[S^2]\otimes_{\bS[S^\lambda]}\bS \right )$ (see \autoref{sec:bokst-peri-spoke})
where $\lambda$ is the representation on
$\mathbb{R}^{2}$ where a chosen generator of $C_p$ rotates by
the angle $2\pi /p$ and $N_e^{C_p}$ is the equivariant
norm constructed in~\cite{HHR16}. The theorem is also valid for $p=2$, where
$\spoke$ should be understood as $\sigma$ and $\mF_{\ell}[S^{1+\sigma}]$ should be understood as the free associative algebra on $S^{1+\sigma}$, recovering \cite[Theorem 5.8]{DMPR21}
(see \autoref{rem:DMPR-thm}).
\autoref{thm:Bok} is proved in \autoref{sec:bokst-peri-spoke}.

\medskip
One of the main difficulties in $C_p$-equivariant homotopy theory is the fact
that the $\RO(C_p)$-graded $C_p$-equivariant dual Steenrod algebra is not a flat Hopf algebroid when $p$ is an odd prime and the algebra structure is complicated~\cite{SW22,HKSZ22}. 
In contrast, we show that the pair $(\pispa\mF_p,
\pispa\sTHH(\mF_p))$  is a flat Hopf algebroid when one considers the $\ROs(C_p)$-graded $a_{\Yright}$-free
  homotopy groups $\pispa$ as developed by \cite{Wil17,HSZZ} and the algebra structure is much simpler. See \autoref{appendixA} for the definition of $\pispa$ and the computation of $\pispa\mF_p$. 
\begin{thm}\label{thm:sTHH}
Let $p$ be an odd prime. We have
\[
(\pispa\mF_p,\pispa\sTHH(\mF_p)=(\bFp[a_{\Yright}, u_{\lambda}]\langle u_{\Yright} \rangle, \bFp[a_{\Yright},
  u_{\lambda}, \mathrm{Nm}(\mu)]\langle u_{\Yright}, \tilde{\mu} \rangle)
\]
where $|\tilde{\mu}|=1+\spoke$ and $|\mathrm{Nm}(\mu)|=2\rho$.
The right unit formula satisfies 
\begin{align}
\eta_R(a_{\Yright})&=a_{\Yright}\\ 
\label{ulambda}  \eta_R(u_{\lambda})& =u_{\lambda}+\beta \cdot  a_{\spoke}^{2p}\mathrm{Nm}(\mu)\\
\label{uspoke}  \eta_R(u_{\spoke})& =  u_{\spoke}+\beta' \cdot  a_{\spoke}^2\tilde{\mu}
\end{align}
where $\beta,\beta'$ are units in $\mathbb{F}_p$.
The coproduct satisfies 
\begin{align*}
  \Delta(\tilde{\mu}) & = \tilde{\mu} \otimes 1 + 1 \otimes \tilde{\mu}\\
  \Delta (\mathrm{Nm}(\mu))& =\mathrm{Nm}(\mu)\otimes 1+1\otimes \mathrm{Nm}(\mu) \,. 
\end{align*}
\end{thm}
\noindent We compute ring structure 
and the Hopf algebroid
structure maps in \autoref{sec:hopf-algebroid}.

\medskip

With the Hopf algebroid from \autoref{thm:sTHH}, we compute the Tate construction of
$N_e^{C_p}\bF_p$, giving a new and simpler proof of the Segal conjecture for $C_p$.
 Our computation uses the descent spectral sequence associated to the map
 $N_e^{C_p}\bF_p\to \mF_p$ of $C_p$-commutative ring spectra, see
 \autoref{descent}, along with a May spectral sequence that computes the
 $\mathrm{E}_2$-page of this descent spectral sequence, see \autoref{mayss}.
 This answers~\cite[Question~6.3]{HW21}, extending work of~\cite{HW21} to arbitrary primes.
\begin{thm}\label{thm:Segal}
The spectrum $N_e^{C_p}\bF_p$ is Borel complete.  
\end{thm}

\begin{cor}[Segal conjecture for $C_p$] 
  For any bounded below spectrum $X$, the Tate diagonal map 
\[ 
X\to (X^{\wedge p})^{tC_p}
\]
is the $p$-completion map. 
\end{cor}

\begin{proof}
As in~\cite[Theorem~III.1.7]{NS18}, the fact that
\[ 
N_e^{C_p}\bF_p\to (N^{C_p}_e\bF_p)_{a_{\spoke}}^{\wedge}
\]
is an equivalence implies that the completion map 
\[ N^{C_p}_eX\to (N^{C_p}_eX)_{a_{\spoke}}^{\wedge}\]
is an equivalence after $p$-completion for a bounded below spectrum $X$.
Taking geometric fixed points and noticing that there are equivalences $(N_e^{C_p}X)^{\Phi C_p} \simeq X$  and $((N^{C_p}_eX)_{a_{\spoke}}^{\wedge})^{\Phi C_p} \simeq (X^{\wedge
  p})^{tC_p}$ gives the claim.
\end{proof}
\noindent This result was originally proven by Gunawardena~\cite{Gun81} for the
sphere spectrum $\mathbb{S}$ and by Lunoe-Nielsen--Rognes~\cite[Theorem~5.13]{LNR12} for general
 $X$ under an additional hypothesis that $X$ has finite type mod $p$
 homology. Nikolaus--Scholze~\cite{NS18} proved
 that the more general statement reduces to \autoref{thm:Segal}.
 
\subsection{Outlook}\label{outlook}
Work in progress of the first author with  Behrens, Johnson and Kong~\cite{AKBJK} produces a spectrum $\mathrm{MUP}_{\mu_p}$ using a $C_p$-equivariant Snaith construction, which has the property that the 
$
\mathrm{MUP}_{\mu_p}^{\Phi C_p}
$
is a $\mathbb{F}_p$-algebra. 
This is closely related to a $C_p$-equivariant analogue of $\mathrm{BP}_{\mathbb{R}}$ constructed by the second author with Hu, Kriz and Somberg~\cite{HKSZ24} and extended to the genuine equivariant setting by~\cite{AKBJK2}. 
Let 
\[
\mathrm{MUP}_{\mu_{p^n}}:=N_{C_p}^{C_{p^n}}\mathrm{MUP}_{\mu_p} \,.
\] 
Then $\mathrm{MUP}_{\mu_{p^2}}^{\Phi C_p}$ is a $N^{C_p}_e\mathbb{F}_p$-algebra. The goal of the program of the first author and collaborators is to approach the $3$-primary Kervaire invariant using using the $C_9$-spectrum $\mathrm{MUP}_{\mu_9}$ or some variation on this spectrum following the suggestion of Hill--Hopkins--Ravenel~\cite{HHRodd}. Understanding the $C_p$-fixed points of $N^{C_p}\mF_p$ is therefore helpful for this goal and our computational results in \autoref{descent} give a first step towards that computation. 

\begin{quest}
Can the descent spectral sequence be used to study $(N_e^{C_3}\mathbb{F}_3)^{C_3}$ and shed light on $\mathrm{MUP}_{\mu_9}$? 
\end{quest}

A $\mathbb{E}_{\spoke}$ algebra is the data of $C_p$-spectrum $R$ such that $i_e^*R$ is an $\mathbb{E}_1$ algebra and $R$ is a $\mathbb{E}_0$ algebra in $N^{C_p}i_e^*R$-bimodules. Although, we only consider $C_p$-commutative ring spectra in this paper, one can more generally consider spoke topological Hochschild homology of a $\mathbb{E}_{\spoke}$ algebra. One might hope to use our theory of spoke THH to study deformations of $\mathbb{E}_{\spoke}$ algebras. Deformation theory of $\mathbb{E}_{\spoke}$ algebras could be helpful for producing orientations along the lines of~\cite{HSW23}.  

\begin{quest}\label{deformation-question}
Can spoke topological Hochschild homology be used to study deformations of $\mathbb{E}_{\spoke}$ algebras? 
\end{quest}

The topological Hochschild homology of spherical group rings was one of the first computations of interest because of the relation to diffeomorphism groups of manifolds~\cite{BHM89}. The interpretation in terms of free loop spaces~\cite{Goo85} also led to applications in string topology. The setting of Real topological Hochschild homology a similar interpretation in terms of signed free loops spaces holds by work of~\cite{Hog16,DMPR21,DMPR24} and in the setting of quaternionic topological Hochschild homology a formula in terms of twisted free loops spaces was proven by the first author, P\'eroux and Merling~\cite{AKMP24}. 

\begin{quest}
Is there a free loop space interpretation of the spoke topological Hochschild homology of spherical group rings?
\end{quest}


\subsection{Conventions}\label{conventions}
 We write 
\begin{itemize}
\item  $p$ and $\ell$ for primes, $C_p$ for the cyclic group of order $p$.
\item $\rho$ for the regular representation of $C_p$.
\item $\lambda$ for a $2$-dimensional irreducible representation of $C_p$.
\item $S^V$ for the one point compactification of a real $C_p$-representation
  $V$.
\item $S^{\spoke}$ for the unreduced suspension of the $C_p$-set $C_p$, which we call the ``spoke sphere''. 
The inclusion of $1$-skeleton fits in a cofiber sequence 
\begin{equation}
\label{eq:cofib}
S^{\Yright} \longrightarrow S^{\lambda} \overset{q}{\longrightarrow} (C_p)_+\otimes S^{2}.
\end{equation}
\item $\Top^{C_p}$ for the category of $C_p$-spaces and $\Sp^{C_p}$ for the category of $C_p$-spectra.
\item $\CAlg_{C_p}$ for the category strict commutative monoids in
  $C_p$-spectra, which is equivalent to algebras over a complete $N_{\infty}$ operad by~\cite[Theorem A.6]{BH15}. Also, by~\cite[Theorem~7.27]{LLP25} this is equivalent to $C_p$-$\mathbb{E}_\infty$-rings in the sense of~\cite{NS22} or normed $G$-algebras in the sense of~\cite{BH21}, but we do not make essential use of this fact. We will call an object $R\in \CAlg_{C_p}$ a $C_p$-$\mathbb{E}_\infty$-ring for brevity in the body of the paper and we called such an object a $C_p$-commutative ring spectrum in the introduction for accessibility. 

\item $\mF_p$ for the (Eilenberg--MacLane spectrum of) the constant Mackey functor of a finite field $\bF_p$ of order $p$. 
\item $\otimes$ for the smash product of ($C_p$-)spectra and tensoring of a ($C_p$-)space with a ($C_p$-)spectrum. We use $\odot$ for the tensoring of a ($C_p$-)-space with a ($C_p$-)$\mathbb{E}_{\infty}$-ring. 
\item $\bF_p[y]$ for a polynomial algebra on an even degree class $y$
  and $\bF_p\langle x\rangle$ for an exterior algebra on an odd degree class $x$. 
\end{itemize}

\subsection{Acknowledgments}
The authors would like to thank Mark Behrens, Eva Belmont, Jeremy Hahn, Michael A. Hill, and Adela (Yiyu) Zhang for helpful conversations and Jens Hornbostel for useful comments on a draft of this paper. In particular, the first author first learned  the definition of spoke topological Hochschild homology from Jeremy Hahn and Dylan Wilson. 
The authors would like to thank the Isaac Newton Institute for Mathematical Sciences, Cambridge, for support and hospitality during the programme \emph{New horizons for equivariance in homotopy theory}, where work on this paper was undertaken. This work was supported by EPSRC grant EP/Z000580/1. 
Angelini-Knoll is grateful to Max Planck Institute for Mathematics in
Bonn for its hospitality and financial support. 

\section{Spoke topological Hochschild homology}
We begin with our key definition. 
\begin{defn}
Let $R$ be a $C_p$-$\bE_{\infty}$-ring in $C_p$-spectra. We define 
\[ 
\sTHH(R):=S^{\spoke}\odot R,
\]
using the tensoring 
\[ \-- \odot\--  : \Top^{C_p}\times \CAlg_{C_p} \longrightarrow  \CAlg_{C_p}
\]
of a $C_p$-space with a $C_p$-$\bE_{\infty}$-ring.
\end{defn}
\begin{rem}
One way to define the tensoring of a $C_p$-space with a $C_p$-$\bE_{\infty}$-ring is to take the $C_p$-colimit of the map 
\[S^{\spoke}\to \mathcal{O}_{C_p}\overset{R}{\longrightarrow} \mathrm{CAlg}_{C_p}\] 
in the sense of~\cite{Sha23} (cf.~\cite[Definition~5.2]{QS21} and \cite[Corollary 3.22]{Ste25}). Alternatively, we can use the equivariant Loday construction of~\cite{LRZLoday}. 
\end{rem}
\begin{rem}
  For $p=2$, we have $S^{\Yright} = S^{\sigma}$ where $\sigma$ is the sign representation and
  $\sTHH = \mathrm{THR}$ defined in \cite{DMPR21} and \cite[Definition~5.2]{QS21}. 
\end{rem}
\begin{rem}
  The restriction is $i_e^{*}(\sTHH(R)) \simeq (\vee_{p-1}S^1) \odot R \simeq \THH(R)^{\otimes_R (p-1)}$.
\end{rem}

We then identify the geometric fixed points of this construction. 
\begin{lem}
There is an equivalence 
\[ \sTHH(R)\simeq R\otimes_{N_e^{C_p}i_e^*R}R\,.\]
Consequently, there is an equivalence 
\[ 
\sTHH(R)^{\Phi C_p}\simeq R^{\Phi C_p}\otimes_{i_e^*R}R^{\Phi C_p} \,.
\]
\end{lem}
\begin{proof}
This is a special case of~\cite[Corollary~3.36]{Ste25}. The last statement follows using the fact that the geometric fixed points functor is symmetric monoidal. 
\end{proof}

The following corollary is immediate. 
\begin{cor}\label{cor:F_p}
There is an isomorphism of $\pi_{*}\mF_p^{\Phi C_p}$-algebras
\[ \pi_{*}\sTHH(\mF_p)^{\Phi C_p}=\bF_p[y,\bar{y}]\langle x,\bar{x}\rangle 
\]
where $|x|=|\bar{x}|=1$ and $|y|=|\bar{y}| = 2$. Here we note that 
$\pi_*\mF_p^{\Phi C_p}=\bF_p[y]\langle x\rangle$ where
 $|y|=2$ and $|x|=1$. The generators are chosen such that the left and right
 unit maps
\[\eta_L\,,\eta_R :\pi_*\mF_p^{\Phi C_p}\longrightarrow \pi_*\sTHH(\mF_p)^{\Phi C_p} \]
induced by the two distinct inclusions of a fixed point $\eta_L, \eta_R: * \to S^{\Yright}$ satisfy $\eta_L(y)=y$, $\eta_L(x)=x$, $\eta_R(x)=\bar{x}$ and $\eta_R(y)=\bar{y}$. Moreover, when $\ell$ is a prime distinct from $p$ then 
$\sTHH(\mF_\ell)^{\Phi C_p}=0$.
\end{cor}

\section{B\"okstedt periodicity for spoke THH}
\label{sec:bokst-peri-spoke}
The goal of this section is to prove~\autoref{thm:Bok}. Here a definition is in order.
For a $C_p$-space $X$, let $\bS[X] = \bS \otimes \left (S^0 \oplus X \oplus X^{\otimes 2} \oplus \cdots
\right )$ be the free $C_p$-equivariant $\mathbb{E}_1$-algebra on $X$.
\begin{defn}[\relax{cf.~\cite[Definition~2.2]{HSW23}}]
  \label{defn:S-Sspoke}
    We define 
\[ \bS[S^{1+\spoke}]:=N_e^{C_p}\bS[S^2]\otimes_{\bS[S^\lambda]}\bS.\]
Here, $N_e^{C_2}\bS[S^2]$ is an $\bE_1$-algebra and a $\bS[S^\lambda]$-bimodule via the $\bE_1$-map induced by the composite
\begin{equation}\label{equiv-map}
S^{\lambda}\overset{q}{\longrightarrow} (C_p)_+\otimes S^{2}\overset{\varphi}{\longrightarrow} N_e^{C_p}\bS[S^2]\,,
\end{equation}
where the first map is the quotient $q$ from \autoref{eq:cofib} and the second
map $\varphi$ is the canonical inclusion to $N_e^{C_p}\bS[S^2] \cong \mathbb{S} \oplus (C_p)_+\otimes S^{2} \oplus \cdots$.
For a $C_p$-$\mathbb{E}_{\infty}$-ring $R$, we write 
\[ R[S^{1+\spoke}]:=R\otimes \bS[S^{1+\spoke}]\,.
\]
Note that $R[S^{1+\spoke}]$ is a unital left module over $R\otimes N_e^{C_p}\bS[S^2]$. 
\end{defn} 

Towards proving the theorem, we first observe that we can understand the underlying homotopy groups equipped with their Weyl group action. 
\begin{lem}
  \label{lem:Weyl}
We can identify 
\[ \pi_{*}^{e}\sTHH(\mF_\ell)=\bF_{\ell}[\mu_{1},\mu_{2},\cdots \mu_{p-1}]\]
where $|\mu_{i}|=2$ for $1\le i\le p-1$ and the Weyl group $C_p$ acts by $\gamma(\mu_i)=\mu_{i+1}$ for $1\le i\le p-2$ and $\gamma(\mu_{p-1})=-\sum_{j=1}^{p-1}\mu_j$. 
\end{lem}
\begin{proof}
We observe that 
\[ 
i_e^*\sTHH(\mF_\ell)\simeq \THH(\bF_\ell)^{\otimes_{\bF_\ell}(p-1)}
\]
so the result, without the $C_p$-action, follows from~\cite{Bok85}. For the $C_p$-action, we label the spokes of $S^\spoke$ by $e_i$ for $0\le i\le p-1$. 
Then the $p-1$ copies of $\THH_*(\bF_\ell)$ correspond to the loops $e_i-e_0$ for $1\le i\le p-1$. Since the generator $\gamma$ of $C_p$ cyclically permutes the edges $e_i$, this determines the action $\gamma (e_i-e_{i-1})=e_{i+1}-e_i$ for 
$1\le i\le p-2$ and $\gamma (e_{p-1}-e_{p-2})=e_{0}-e_{p-1}$.
The generator $\mu_i$ corresponds to $e_i-e_{i-1}$, and
$e_0-e_{p-1}=e_0-e_1+e_1-e_2+ \ldots +e_{p-2}-e_{p-1}$
corresponds to $-\sum_{j=1}^{p-1}\mu_j $
as desired.
\end{proof}

The following construction is a special case of a construction from~\cite[\S~6]{AKKQ25}, which is an equivariant analogue of~\cite[\S~A.1]{HW22}. 
\begin{con}[The $(1+\spoke)$-suspension map]\label{suspension-map}
Let $R$ be a $C_p$-$\bE_{\infty}$-ring. Consider the $C_p$-space of non-equivariant maps $*\to S^{\spoke}$. Each
  map induces a map $R \to S^{\Yright} \odot R$ in the $C_p$-space of non-equivariant maps
  of spectra, and thus there is a map

  \[
   S^{\Yright} \simeq \mathrm{Map}(*,S^{\spoke})\to \mathrm{Map}(R,\sTHH(R)) 
\]
with adjoint map 
\[ S^{\spoke}_+\otimes R\to \sTHH(R)\,.\]
This is natural in $R$, producing a diagram 
\[
\begin{tikzcd}
S^{\spoke}_+\otimes \bS \simeq (S^{\Yright}\otimes \bS) \oplus \bS \arrow{r} \arrow{d} &  \sTHH(\bS)\simeq\bS  \arrow{d} \\ 
S^{\spoke}_+\otimes R\simeq (S^{\Yright}\otimes R) \oplus R \arrow{r} & \sTHH(R)
\end{tikzcd}
\]
and using the canonical splitting of the top map we produce a commutative diagram
\[
\begin{tikzcd}
S^{\spoke}\otimes \bS \arrow{r} \arrow{d} &  0 \arrow{d} \\ 
S^{\spoke}\otimes R \arrow{r} & \sTHH(R)
\end{tikzcd}
\]
and consequently a map 
\[ 
\sigma^{1+\spoke} :S^{1+\spoke} \otimes \mathrm{fib}(\bS\to R)=S^{\spoke} \otimes \mathrm{cof}(\bS\to R) \longrightarrow \sTHH(R) 
\]
which we call the $(1+\spoke)$-suspension map. 
\end{con}

\begin{rem}\label{suspension-restriction}
The construction above can also be done in the setting of topological spaces and $\bE_\infty$-rings in spectra producing a map 
\[ \bigoplus_{p-1}\mathbb{S}^1\otimes \mathrm{cof}(\mathbb{S}\to R)\longrightarrow \THH(R)^{\otimes_{R} p-1 } 
\]
which is the reduced suspension map of~\cite[Example~A.2.4]{HW22} in each coordinate. 
\end{rem}

\begin{lem}\label{mus}
Let $p$ and $\ell$ be primes. We can choose linearly independent generators 
\[
\mu_1,\cdots ,\mu_{p-1}\in \pi_2^e\sTHH(\mF_\ell)\] such that there is  $C_p$-equivariant refinement 
\[
\tilde{\mu}: S^{1+\spoke}\longrightarrow \sTHH(\mF_\ell)
\]
and this lifts to a map
\begin{equation}
\label{eq:compare-sTHH}
\mF_\ell[S^{1+\spoke}]\longrightarrow \sTHH(\mF_\ell)\,.
\end{equation}
\end{lem}

\begin{proof}
We take the $(1+\spoke)$-suspension map 
\[ 
\sigma^{1+\spoke} :S^{1+\spoke}\otimes \mathrm{fib}(\bS\to \mF_\ell)\longrightarrow \sTHH(\mF_\ell)
\]
from~\autoref{suspension-map}. Since $\ell\in \pi_0^{C_p}\bS=\bZ[t]/(t^2-pt)$ maps to zero in $\pi_0^{C_p}\mF_\ell=\bF_\ell$, there is a non-trivial class 
\[ \ell\in \pi_0^{C_p}\textup{fib}(\bS\to\mF_\ell).\]
The image of it is a class
\[ 
\sigma^{1+\spoke}(\ell)\in \pi_{1+\spoke}^{C_p}\sTHH(\mF_\ell)= \pi_{0}^{C_p}\mathrm{Map}(S^{1+\spoke},\sTHH(\mF_\ell)).
\]
This class is represented by a $C_p$-equivariant map 
\[
  \tilde{\mu} : S^{1+\spoke}\longrightarrow \sTHH(\mF_\ell).
\]
The restriction  
\[
(\tilde{\mu})^e: (S^{1+\Yright})^e \simeq \bigvee_{p-1}S^{2}\longrightarrow \sTHH(\mF_\ell)^e 
\]
picks out some choice of linearly independent classes in $\pi_2^e\sTHH(\mF_\ell)$ by~\autoref{suspension-restriction}. Adjusting our choice of generators, we can choose $\mu_1,\cdots,\mu_{p-1}$ to be given by this choice of linearly independent classes (in particular, these can be chosen to be the same as our previous choice of classes up to a uniform choice of unit). 

We then construct the map \autoref{eq:compare-sTHH}. 
We can rotate the cofiber sequence from \autoref{eq:cofib} to produce the cofiber sequence:
\begin{equation}\label{eq:cofib-rot}
  \begin{tikzcd}
    S^{\lambda} \ar[r,"q"] & (C_p)_+ \otimes S^2 \ar[r,"\partial"] & S^{1+\Yright}
  \end{tikzcd}
\end{equation}
We include  
$\iota: S^{2}\overset{}{\longrightarrow} ((C_p)_+\otimes S^2)^e$ via inclusion of the neutral element into
$C_p$. This yields a map
\begin{equation}
\label{eq:5}
\bS[S^{2}] \overset{\bS[\iota]}{\longrightarrow} \bS[((C_p)_+\otimes S^2)^e] \overset{\bS[\partial^e]}{\longrightarrow} \bS[(S^{1+\Yright} )^e] \overset{\bar{\mu}^{e}}{ \longrightarrow }\sTHH(\mF_\ell)^e,
\end{equation}
where  $\bar{\mu}$ is the extension of $\tilde{\mu}$ using the $\bE_1$-structure of $\sTHH(\mF_\ell)$.
Applying the norm to it and using the fact that $\THH^{\spoke}(\bF_\ell)$ is a $C_p$-$\bE_{\infty}$-ring, we produce 
\[ 
\xi: N_e^{C_p}\bS[S^2] \longrightarrow N_e^{C_p}i_e^*\sTHH(\mF_\ell) \longrightarrow\sTHH(\mF_\ell)\,.
\]
We claim that $\xi$ composed with \autoref{equiv-map} is null.  If so, we obtain a map
\[ 
\bS[S^{1+\spoke}]\longrightarrow \sTHH(\mF_\ell).
\]
Since the target is a $\mF_\ell$-algebra, we can extend scalars to produce the
map \autoref{eq:compare-sTHH}.

Now we prove the claim. We can apply norm to each step of \autoref{eq:5} and
factor the map $\xi$ as below:
\begin{equation*}
\begin{tikzcd}
    S^{\lambda} \ar[r,"q"] \ar[d,"\autoref{equiv-map}"' ] & (C_p)_+ \otimes S^2 \ar[r, "\Delta\otimes\mathrm{id}"]  \ar[ld,"\varphi"] 
    &(C_p)_+ \otimes (C_p)_+ \otimes S^2   \ar[r," \mathrm{id}\otimes\partial"] \ar[ld,"\varphi"]  & (C_p)_+ \otimes
    S^{1+\Yright}\ar[ld,"\varphi"] \ar[d, "\triangledown"]\\
    N_e^{C_p}\bS[S^2] \ar[r] \ar[rrd, out=-30, in=175, "\xi"] & N_e^{C_p}\bS[(C_p)_+ \otimes S^2] \ar[r] &
    N_e^{C_p}\bS[(S^{1+\Yright})^e] \ar[d, "\zeta"] \ar[r, phantom, "\textup{\textcircled{x}}"]& S^{1+\Yright} \ar[ld, "\tilde{\mu}"] \\
     & &  \sTHH(\mF_\ell)
\end{tikzcd}
\end{equation*}
The maps $\varphi$ are inclusions of cells as described below
\autoref{equiv-map}, and $\triangledown$ is adjoint to the identity map. As $\zeta$ extends
$\tilde{\mu}$, the diagram \textcircled{x} commutes. The composite from $S^{\lambda}$
to $S^{1+\Yright}$ fits in the diagram below
\begin{equation*}
  \begin{tikzcd}
    S^{\lambda} \ar[r,"q"]  & (C_p)_+ \otimes S^2 \ar[r, "\Delta\otimes\mathrm{id}"]  \ar[rd, dotted, " \mathrm{id}"']
    &(C_p)_+ \otimes (C_p)_+ \otimes S^2   \ar[r," \mathrm{id}\otimes \partial"]  \ar[d, "\triangledown"] & (C_p)_+ \otimes
    S^{1+\Yright}\ar[d, "\triangledown"]\\
    & & (C_p)_+ \otimes S^2 \ar[r,"\partial"] & S^{1+\Yright} 
  \end{tikzcd}
\end{equation*}
and the dotted composite is the identity
map. Therefore the composite map is null from the cofiber sequence \autoref{eq:cofib-rot}, which
finishes the proof.
\end{proof}

 In the case $p=\ell$, we describe preferred generators for $\mF_p[S^{1+\spoke}]$.  
\begin{con}[{cf.~\cite[Construction~4.7]{HSW23}}]
  \label{const:norm_mu}
 Consider the commutative diagram, where the top row is a cofiber sequence \autoref{eq:cofib-rot}; the top vertical maps are
inclusions of cells; the bottom maps are from \autoref{defn:S-Sspoke}:
\begin{equation*}
    \begin{tikzcd}
        S^{\lambda}\ar[r,"q"] \ar[d,"x"'] & {C_p}_+\otimes S^2 \ar[d] \ar[r,"\partial"] &
        S^{1+\spoke} \ar[ldd,dotted, bend left, "\tilde{x}"] \\
        \bS[S^\lambda] \ar[r] \ar[d] & N_e^{C_p}\bS[S^2] \ar[d] \\
        \bS \ar[r] & \bS[S^{1+\spoke}]
    \end{tikzcd}
\end{equation*}
Since the left vertical composite is null, we produce the dotted map from the cofiber
\[ \tilde{x} : S^{1+\spoke}\to \bS[S^{1+\spoke}]\,.
\]
On the other hand, applying the norm $N_e^{C_p}i_e^*$ to $x$ produces a map 
\[ \mathrm{Nm}(x) : S^{2\rho}\to N_e^{C_p} \bS[S^2].
\]
Finally, using the module structure of $\mathrm{N}_e^{C_p}\bS[S^2]$ on
$\bS[S^{1+\spoke}]$, we produce maps
\[ 
\mathrm{Nm}(x)^k\tilde{x}^{\varepsilon} : S^{2\rho k +\varepsilon (1+\spoke)} \longrightarrow \bS[S^{1+\spoke}]
\]
for $k\ge 0$ and $\varepsilon\in \{0,1\}$. 
\end{con}

\begin{proof}[Proof of~\autoref{thm:Bok}]
We want to show that the map 
\[ 
\mF_\ell[S^{1+\spoke}]\longrightarrow \sTHH(\mF_\ell)\,.
\]
from~\autoref{mus} is a $C_p$-equivariant equivalence. Since the pair $((\--)^e,(\--)^{\Phi C_p})$ of functors is conservative, it suffices to prove an equivalence on underlying and geometric fixed points. On underlying, we observe that 
\[ 
\pi_*^e(\mF_\ell[S^{1+\spoke}])\cong \bF_\ell[m_1,\cdots,m_p] \otimes_{\bF_\ell[m]}\bF_\ell \cong \bF_\ell[m_1,m_2,\cdots ,m_{p-1}] 
\]
where $|m|=2$, $|m_i|=2$ for $1 \leq i \leq p$ and $m \mapsto m_1 + \cdots + m_p$ (and the Weyl group action is the same as in \autoref{lem:Weyl}).
By construction, the map 
\[ 
\pi_*^e\mF_\ell[S^{1+\spoke}]\longrightarrow \pi_{*}^e\sTHH(\mF_\ell)
\]
sends $m_i$ to $\mu_i$. 

When $\ell\ne p$, the source and target are both zero on geometric fixed points so the map is an equivalence also on geometric fixed points. This proves the result when $\ell\ne p$. 

Suppose $\ell=p$. On geometric fixed points, we determine that 
\begin{equation}
\label{eq:geo-of-Sspoke}
 \pi_*\left (\mF_p[S^{1+\spoke}])^{\Phi C_p}\cong \pi_*(\mF_p^{\Phi C_p} \otimes \left (\bS[S^2]\otimes_{\bS[S^0]}\bS \right )\right )\cong \bF_p[u,v]\langle a,b\rangle 
\end{equation}
where $|u|=|v|=2$ and $|a|=|b|=1$. By~\autoref{cor:F_p}, this is abstractly isomorphic to $\sTHH(\mF_p)^{\Phi C_p}$ so it suffices to demonstrate that the map from~\autoref{mus} induces an isomorphism.  

We first claim that the classes $u$ and $a$ map to $y$ and $x$ respectively. To see this, note that that $\mF_p$ splits off of $\mF_p\otimes \bS[S^{1+\Yright}]$ and $\sTHH(\mF_p)$ compatibly and we can choose our classes $u$, $a$, $y$ and $x$ to be compatible with this splitting in the sense that the composite 
\[ 
\pi_*\mF_p^{\Phi C_p}\to \pi_{*}(\mF_p\otimes \bS[S^{1+\Yright}])^{\Phi C_p} \longrightarrow \pi_{*}\sTHH(\mF_p)^{\Phi C_p}\to \pi_*\mF_p^{\Phi C_p}
\]
is the identity map.

The classes $v$ and $a$ denote the Hurewicz image of classes
$v\in\pi_2(\bS[S^2]\otimes_{\bS[S^0]}\bS)$ and $b\in \pi_1(\bS[S^2]\otimes_{\bS[S^0]}\bS)$. By \autoref{const:norm_mu}, these classes can be described
as $\mathrm{Nm}(\mu)^{\Phi C_p}$ and $\tilde{\mu}^{\Phi C_p}$ respectively. By
construction of the map from \autoref{mus}, these classes have non-trivial image
in  $\sTHH(\mF_p)^{\Phi C_p}$. 
We then observe that the image of $\mathrm{Nm}(\mu)^{\Phi C_p}$ and $\tilde{\mu}^{\Phi C_p}$ are trivial in the image of the map
\[
\sTHH(\mF_p)^{\Phi C_p}\to \mF_p^{\Phi C_p}
\] 
and therefore they must correspond to $\overline{y}$ and $\overline{x}$ respectively up to a change of basis, in other words the map is an isomorphism. 
\end{proof}

\begin{rem} \label{rem:DMPR-thm}
 For $p=2$, $\Yright$ shall be understood as $\sigma$ and $\lambda$ as $2\sigma$. There are a prior two definitions for
 $\bS[S^{1+\spoke}]$: one is the free $\bE_1$-algebra on $S^{1+\sigma}$, which we
 will denote as $\bS[S^{\rho}]$; the other is
 $N_e^{C_2}\bS[S^2]\otimes_{\bS[S^{2\sigma}]}\bS$, which we will denote as $\bS[S^{1+\Yright}]$.
 All constructions and propositions above are also valid,
 with the only changes in the proof being that
 $\pi_*\mF_2^{\Phi C_2} = \bF_2[x]$ where $|x|=1$ and that
 \autoref{eq:geo-of-Sspoke} becomes 
$$ \pi_*\left (\mF_2[S^{1+\spoke}] \right)^{\Phi C_2}\cong \bF_2[a,v]\langle b\rangle $$
where $|a|=1$, $|v|=2$ and $|b|=1$. Note that this is not an isomorphism of
rings, but of $\bF_2[a,v]$-modules. We still have that
$\pi_*\left (\mF_2[S^{1+\spoke}] \right )^{\Phi C_2} \to \pi_{*}\THR(\mF_2)^{\Phi  C_2}$ is an
isomorphism. Together with the first statement of \autoref{prop:THH-spoke} below, which is
also valid for $p=2$, we obtain 
\begin{equation*}
\mF_\ell[S^{\rho}] \simeq \mF_\ell[S^{1+\Yright}] \simeq \THR(\mF_\ell).
\end{equation*}

\end{rem}

\section{The Hopf algebroid}
\label{sec:hopf-algebroid}
In this section, we fix an odd prime and determine some formulas for the Hopf algebroid 
\begin{equation}\label{Hopf}(\pispa\mF_p,\pispa\sTHH(\mF_p))
\end{equation}
as claimed in \autoref{thm:sTHH}. See \autoref{appendixA} for an explanation of
the notations $a_{\Yright}$, $u_{\Yright}$, $u_{\lambda}$, $\pisp$ and $\pispa$ as well as a proof of $\pispa(\mF_p) = \bFp[a_{\Yright},  u_{\lambda}]\langle u_{\Yright} \rangle$.
\begin{prop} \label{prop:THH-spoke}
There is an equivalence 
\[\mF_p \otimes \left ( \bigoplus_{\substack{k\ge 0\\ \varepsilon\in \{0,1\}}}\bS^{2\rho k +\varepsilon
      (1+\spoke)} \right )
  \oplus \mF_p \otimes \left( \bigoplus_{\substack{k \ge 2 \\ 1 \leq i \leq m_k}} \Sigma^{2k} (C_p)_+\right) \overset{\simeq}{\longrightarrow}  \sTHH(\mF_p)  \,,
\]
where $m_k = \lfloor \frac{\binom{k+p-2}{p-2}}{p} \rfloor$.  

Consequently, we can identify 
\[  
\pispa\sTHH(\mF_p)=\bFp[a_{\Yright},  u_{\lambda}, \mathrm{Nm}(\mu)]\langle u_{\Yright}, \tilde{\mu} \rangle \,.
\]
where $|\tilde{\mu}|=1+\spoke$ and $|\mathrm{Nm}(\mu)|=2\rho$
as $\pispa\mF_p$-algebras. 
\end{prop}
\begin{proof}
The proof of~\cite[Theorem 4.6]{HSW23} carries over producing an equivalence
\[\bigoplus_{\substack{k\ge 0\\ \varepsilon\in \{0,1\}}}\bS^{2\rho k +\varepsilon (1+\spoke)}
   \oplus \bigoplus_{\substack{k \ge 2 \\ 1 \leq i \leq m_k}} \Sigma^{2k} (C_p)_+ \overset{\simeq}{\longrightarrow} \bS[S^{1+\spoke}]
\]
on $C_p$-geometric fixed points.
Here the number $m_k$ is the number of regular representations in
  $H_{2k}^e(\bS[S^{1+\spoke}], \bFp) \cong \mathrm{Sym}^k_{\bFp}(\bar{\rho})$
  \cite[III 3.4 - 3.6]{AlmkvistFossum} and~\cite[Proposition~4.10]{HSW23}.
Since we are only considering the $a$-free part, the second statement follows from~\autoref{thm:Bok}. It suffices to show that
$\tilde{\mu}^2=0$. Note that $|\tilde{\mu}^2|=2+\lambda$ and the only non-trivial class
in this degree up to units is $a_{\lambda}^{p-2}\mathrm{Nm}(\mu)$. However, we showed in the proof of \autoref{thm:Bok} that $(\tilde{\mu}^{\Phi C_p})^{2}=0$ whereas $(a_{\lambda}^{p-2}\mathrm{Nm}(\mu))^{\Phi C_p}=(\mathrm{Nm}(\mu))^{\Phi C_p}\ne 0$. 
\end{proof}

We can now prove the desired formulas for the Hopf algebroid \eqref{Hopf}. 
\begin{proof}[Proof of \autoref{thm:sTHH}]
We follow the strategy from~\cite[Theorem~2.3]{HW21}. 
When $R=\bS$, then
\[ (\pispa\mathbb{S},\pispa\sTHH(\bS)) =(\pispa\mathbb{S},\pispa\mathbb{S})
\]
is the trivial Hopf algebroid and therefore $\eta_L=\eta_R=\textup{id}$. Since
$a_{\lambda}$ is in the image of the Hurewicz map 
\[
\pispa\bS\to
\pispa\mF_p\,,
\] 
we have $\eta_R(a_{\lambda}) = a_\lambda$ and therefore
$\eta_R(a_{\Yright}) = a_{\Yright}$. 

For the classes $u_{\lambda}$ and $u_{\Yright}$, we note that 
\begin{align*}
  \eta_R(u_{\lambda}) & \in  \pi_{2-\lambda,\textup{$a$-free}}^{C_p}\sTHH(\mF_p)=\bF_p\{u_{\lambda}, a_{\spoke}^{2p}\mathrm{Nm}(\mu)\}\\
\eta_R(u_{\Yright}) & \in  \pi_{1-\Yright,\textup{$a$-free}}^{C_p}\sTHH(\mF_p)=\bF_p\{u_{\Yright}, a_{\spoke}^2 \tilde{\mu}\}
\end{align*}
so
\begin{align*}
  \eta_R(u_{\lambda}) & =\alpha u_{\lambda}+\beta a_{\spoke}^{2p}\mathrm{Nm}(\mu) \\
  \eta_R(u_{\Yright}) & =\alpha' u_{\Yright}+\beta' a_{\spoke}^2 \tilde{\mu}
\end{align*}
for some $\alpha,\beta, \alpha', \beta' \in \bF_p$.  Since $\mathrm{res}_e^{C_p}(u_{\lambda})=1$ and
 $\mathrm{res}_e^{C_p}(a_{\Yright})=0$ we determine that $\alpha=1$. 
The underlying map of $u_{\Yright}$ is non-trivial by \autoref{lem:u_spoke}. As
the underlying map of $\eta_R: \Ff \to \sTHH(\mF_p) = S^{\Yright} \odot \Ff$ is
homotopic to the underlying map of $\eta_L$, we also have that $\alpha'=1$.

On geometric fixed points, the map $N_e^{C_p}i_e^*\Ff\to \Ff$ is the map
$\Fp\to \Fp[y]\langle x\rangle$ where $y=u_{\lambda}/a_{\lambda}$ and $x=u_{\Yright}/a_{\Yright}$.
We can consider the associated descent Hopf algebroid 
\[ (\bF_p[y]\langle x\rangle , \bF_p[y]\langle x\rangle \otimes \bF_p[y]\langle x\rangle ).
\]
Writing $y=y\otimes1$, $\bar{y} = 1\otimes y$ and writing $x=x\otimes 1$, $\bar{x}=1\otimes x$,
we compute $\eta_R(y) - \eta_L(y) = \bar{y} - y \neq 0$ and $\eta_R(x) - \eta_L(x)= \bar{x}
-x \neq 0$. So $\beta \neq 0$ and $\beta' \neq 0$.

The element $\mathrm{Nm}(\mu)$ is primitive for degree reasons since the only
elements in $\pi_{2\rho,\textup{$a$-free}}^{C_p}(\sTHH(\bF_p)\otimes_{\bF_p}\sTHH(\bF_p))$ are $\mathrm{Nm}(\mu)\otimes
1$ and $1\otimes \mathrm{Nm}(\mu)$.
The element $\tilde{\mu}$ is also primitive for degree reasons. 
\end{proof}

\section{The descent spectral sequence}\label{descent}
We now consider the descent spectral sequence for the map 
\[ 
N^{C_p}_e\bF_p\to \mF_p 
\]
of $C_p$-$\bE_{\infty}$-rings. For a computable $\mathrm{E}_2$-page, we will
  consider it after $a_{\Yright}$-completion and inverting $a_{\spoke}$.
  
\begin{prop}\label{prop:ss}
There is a multiplicative spectral sequence 
\begin{align}\label{SSa-inverted} \tag{DSS}
\mathrm{E}_2^{i+k\spoke,s}\implies \pi_{i+k\spoke}^{C_p}(N_e^{C_p}\mathbb{F}_p)_{a_{\spoke}}^{\wedge}[a_{\spoke}^{-1}]
\end{align}
where 
\[
\mathrm{E}_2^{i+k\spoke,s}=\left (\lim_n \mathrm{Ext}_{\bF_p[\mathrm{Nm}(\mu)]\langle \tilde{\mu} \rangle/(\mathrm{Nm}(\mu)^{p^n})}^{s,s+i+k\spoke}(\bF_p,\bF_p[a_{\spoke},
u_{\lambda}^{\pm 1}]\langle u_{\spoke}\rangle) \right )[a_{\spoke}^{-1}]\,.
\]
Here $\bF_p[\mathrm{Nm}(\mu)]\langle \tilde{\mu} \rangle/(\mathrm{Nm}(\mu)^{p^n})$ is a
primitively generated Hopf algebra over $\bF_p$ and the comodule structure of $\bF_p[a_{\spoke},
u_{\lambda}^{\pm 1}]\langle u_{\spoke}\rangle$ is determined by \autoref{thm:sTHH}.
\end{prop}

\begin{proof}
There is a $\mathbb{E}_{\infty}$-algebra in filtered spectra 
\[ \left ( \lim_{\Delta} \tau_{\ge *} (\underline{\mathbb{F}}_p^{\otimes_{N_e^{C_p}\mathbb{F}_p} \bullet+1}) \right )_{a_{\spoke}}^{\wedge } [a_{\spoke}^{-1}]\,, 
\]
and we consider the associated spectral sequence. 

We first discuss convergence. 
Note that 
\begin{multline}
 \textup{colim}_ n \left ( \left ( \lim_{\Delta} \tau_{\ge -n} (\underline{\mathbb{F}}_p^{\otimes_{N_e^{C_p}\mathbb{F}_p} \bullet+1})\right )_{a_{\spoke}}^{\wedge }  [a_{\spoke}^{-1}] \right ) \\
 \simeq 
 \left ( \textup{colim}_ n \left ( \left ( \lim_{\Delta} \tau_{\ge -n}
  (\underline{\mathbb{F}}_p^{\otimes_{N_e^{C_p}\mathbb{F}_p} \bullet+1}) \right ) _{a_{\spoke}}^{\wedge } \right ) \right )[a_{\spoke}^{-1}]
  \\
\simeq \left (  \left (\lim_{\Delta}  (\underline{\mathbb{F}}_p^{\otimes_{N_e^{C_p}\mathbb{F}_p} \bullet+1}) \right )_{a_{\spoke}}^{\wedge } \right ) [a_{\spoke}^{-1}] 
\simeq (N_e^{C_p}\mathbb{F}_p)_{a_{\spoke}}^{\wedge}[a_{\spoke}^{-1}]
\end{multline}
where the first equivalence follows because colimits commute with colimits, the second equivalence follows by the same argument as \cite[Lemma~I.2.6(ii)]{NS18}. More explicitly, note that the cofiber of the map 
\[ 
 \left (\lim_{\Delta} \tau_{\ge -n} (\underline{\mathbb{F}}_p^{\otimes_{N_e^{C_p}\mathbb{F}_p} \bullet+1}) \right )_{a_{\spoke}}^{\wedge } \longrightarrow \left ( \lim_{\Delta} (\underline{\mathbb{F}}_p^{\otimes_{N_e^{C_p}\mathbb{F}_p} \bullet+1}) \right )_{a_{\spoke}}^{\wedge } 
\]
is 
\[  \left (\lim_{\Delta} \tau_{\le -n-1} (\underline{\mathbb{F}}_p^{\otimes_{N_e^{C_p}\mathbb{F}_p} \bullet+1}) \right )_{a_{\spoke}}^{\wedge}
\]
which has homotopy groups bounded above by $-n$ and therefore it is zero after applying $\colim_n$. 
The last equivalence follows as in~\cite[
Proposition~3]{GIKR22} using~\cite[Theorem~6.6] {Bou79} and its generalization in~\cite{Man24}. 

To see that the limit satisfies 
\[ \lim_n \left( \left (\lim_{\Delta} \tau_{\ge -n}
      (\underline{\mathbb{F}}_p^{\otimes_{N_e^{C_p}\mathbb{F}_p} \bullet+1}) \right
    )_{a_{\spoke}}^{\wedge }[a_{\spoke}^{-1}] \right) \simeq 0,
\]
we first note that 
\[ \lim_n \left (\lim_{\Delta} \tau_{\ge -n}
    (\underline{\mathbb{F}}_p^{\otimes_{N_e^{C_p}\mathbb{F}_p} \bullet+1}) \right
  )_{a_{\spoke}}^{\wedge } \simeq 0
\]
since $a_{\spoke}$-completion commutes with limits. The fiber of
\[ \lim_n \left (\lim_{\Delta} \tau_{\ge -n}
    (\underline{\mathbb{F}}_p^{\otimes_{N_e^{C_p}\mathbb{F}_p} \bullet+1}) \right
  )_{a_{\spoke}}^{\wedge } \longrightarrow \lim_n \left( \left (\lim_{\Delta} \tau_{\ge -n}
      (\underline{\mathbb{F}}_p^{\otimes_{N_e^{C_p}\mathbb{F}_p} \bullet+1}) \right
    )_{a_{\spoke}}^{\wedge }[a_{\spoke}^{-1}] \right) 
\]
can be identified with 
\[ \lim_n {EC_p}_+\otimes \left (\lim_{\Delta} \tau_{\ge -n} (\underline{\mathbb{F}}_p^{\otimes_{N_e^{C_p}\mathbb{F}_p} \bullet+1}) \right )
\]
which is also zero since the connectivity of 

\[
{EC_p}_+\otimes \left (\lim_{\Delta} \tau_{\ge -n} (\underline{\mathbb{F}}_p^{\otimes_{N_e^{C_p}\mathbb{F}_p} \bullet+1}) \right )
\]
is increasing as $n$ goes to $-\infty$. 

We now prove the identification of the $\mathrm{E}_2$-page, which is the
associated graded of the filtration above. We first note that 
\begin{equation*}
\underline{\mathbb{F}}_p^{\otimes_{N_e^{C_p}\mathbb{F}_p} \bullet+1}
 \simeq \sTHH(\mF_p)^{\otimes_{\mF_p} \bullet} .
\end{equation*}
By \autoref{thm:sTHH} we have that $\pisp(\mF_p)_{a_{\spoke}}^{\wedge}[a_{\spoke}^{-1}]=\bF_p[a_{\spoke}^{\pm 1},u_{\lambda}^{\pm 1}]\langle u_{\spoke} \rangle$, and we can identify
\[ 
\pisp\left ( \sTHH(\mF_p)^{\otimes_{\mF_p}q}  \right )_{a_{\spoke}}^{\wedge}[a_{\spoke}^{-1}]
\]
with 
\[
\left ( \bF_p[a_{\spoke},u_{\lambda}^{\pm 1},\mathrm{Nm}(\mu)_1,\cdots ,\mathrm{Nm}(\mu)_q]\langle u_{\spoke}, \tilde{\mu}_1,\cdots, \tilde{\mu}_{q}\rangle \right )_{a_{\spoke}}^{\wedge}[a_{\spoke}^{-1}].
\]

By the postponed \autoref{lem:cobar}, the cobar complex
\[
\pisp \left ( \sTHH(\mF_p)^{\otimes_{\mF_p}\bullet} \right )_{a_{\Yright}}^{\wedge}[a_{\spoke}^{-1}]
\]  
is the limit of cobar complexes 
\[ 
(\lim_n C^{\bullet}_{\bF_p[\mathrm{Nm}(\mu)]\langle \tilde{\mu}\rangle/(\mathrm{Nm}(\mu)^{p^n})}(\bF_p[a_{\spoke},
u_{\lambda}^{\pm 1}]\langle u_{\spoke}\rangle))[a_{\spoke}^{-1}],
\]
and since these are finite abelian groups in each tri-degree, the relevant $R^1\lim$
term vanishes, producing the desired cohomology as the $\mathrm{E}_2$-term of the
  spectral sequence.
\end{proof}

\begin{lem}\label{lem:cobar}
The term
\begin{equation}
\label{eq:3}
 \left ( \bF_p[a_{\spoke},u_{\lambda}^{\pm 1},\mathrm{Nm}(\mu)_1,\cdots ,\mathrm{Nm}(\mu)_q]\langle u_{\spoke}, \tilde{\mu}_1,\cdots, \tilde{\mu}_{q}\rangle \right )_{a_{\spoke}}^{\wedge}
\end{equation}
 can be identified with 
\begin{equation}
\label{eq:4}
\left (\lim_n \bF_p[a_{\spoke},u_{\lambda}^{\pm 1}]\langle u_{\spoke}\rangle\otimes  \bF_p[\mathrm{Nm}(\mu)]\langle \tilde{\mu}\rangle/(\mathrm{Nm}(\mu)^{p^n})^{\otimes q} \right )\,.
\end{equation}
\end{lem}

\begin{proof}
Let \[ J=(a_{\spoke}^8u_{\lambda}^{-1}\mathrm{Nm}(\mu)_1,\cdots ,a_{\spoke}^8u_{\lambda}^{-1}\mathrm{Nm}(\mu)_q, a_{\spoke}^2u_{\lambda}^{-1}u_{\spoke}\tilde{\mu}_1,\cdots ,a_{\spoke}^2u_{\lambda}^{-1}u_{\spoke}\tilde{\mu}_q)\]
be the degree 0 part of the ideal $(a_{\Yright})$.
Since the completion is taken in the graded sense, \eqref{eq:3} can be identified with 
\[ 
\bF_p[a_{\spoke},u_{\lambda}^{\pm 1},\mathrm{Nm}(\mu)_1,\cdots ,\mathrm{Nm}(\mu)_q ]\langle u_{\spoke}, \tilde{\mu}_1,\cdots ,\tilde{\mu}_q\rangle _{J}^{\wedge}.
\]
Similarly, this can be further identified with 
\begin{equation}\label{eq:2}
\bF_p[a_{\spoke},u_{\lambda}^{\pm 1},\mathrm{Nm}(\mu)_1,\cdots,\mathrm{Nm}(\mu)_q]\langle u_{\spoke},\tilde{\mu}_1,\cdots ,\tilde{\mu}_q\rangle_{K}^{\wedge}
\end{equation} 
 where 
 \[ K=(\mathrm{Nm}(\mu)_1,\cdots ,\mathrm{Nm}(\mu)_{q}, \tilde{\mu}_1,\cdots ,\tilde{\mu}_q ),
 \]
as $J$ is also the degree 0 part of $K$. We claim that \eqref{eq:2} can be identified with \eqref{eq:4}. Explicitly, setting $A_n :=\bF_p[a_{\spoke},u_{\lambda}^{\pm 1}]\langle u_{\spoke}\rangle\otimes  \bF_p[\mathrm{Nm}(\mu)]\langle \tilde{\mu}\rangle/(\mathrm{Nm}(\mu)^{p^n})^{\otimes q} $ which can be identified with 
\[
    \bF_p[a_{\spoke},u_{\lambda}^{\pm 1}]\langle u_{\spoke}\rangle\otimes  \bF_p[\mathrm{Nm}(\mu)_1,\cdots
  ,\mathrm{Nm}(\mu)_q]\langle \tilde{\mu}_1, \cdots, \tilde{\mu}_q\rangle/(\mathrm{Nm}(\mu)_1^{p^n},\cdots ,\mathrm{Nm}(\mu)_q^{p^n})
\]
and 
\[
B_n=\bF_p[a_{\spoke},u_{\lambda}^{\pm 1},\mathrm{Nm}(\mu)_1,\cdots,\mathrm{Nm}(\mu)_q]\langle u_{\spoke},\tilde{\mu}_1,\cdots ,\tilde{\mu}_q\rangle/K^{p^n}
\]
then it is clear that 
\[ \lim_n B_n= \bF_p[a_{\spoke},u_{\lambda}^{\pm
    1},\mathrm{Nm}(\mu)_1,\cdots,\mathrm{Nm}(\mu)_q]\langle u_{\spoke}, \tilde{\mu}_1,\cdots ,\tilde{\mu}_q\rangle_{K}^{\wedge}\]
and it suffices to prove that there is a level map 
\[ \{ A_n\} \to \{ B_n\} \]
of pro-abelian groups that is a pro-isomorphism. By~\cite[Lemma~2.3]{Isa01}, it suffices to note that for any integer $s$ there exists an integer $t$ such that 
\[ 
\begin{tikzcd}
A_s \ar[r] \ar[d] & B_s \ar[d]  \ar[dl] \\ 
A_t \ar[r] & B_t  \,.    
\end{tikzcd}
\]
To see this note that there is a level map induced by the inclusion 
\[
(\mathrm{Nm}(\mu)_1^{p^n},\cdots ,\mathrm{Nm}(\mu)_q^{p^n}) \subset K^{p^n}
\]
and that for any integer $s$ there exists an integer $t$ (for example take $t=s+1$) such that 
\[ K^{p^t}\subset (\mathrm{Nm}(\mu)_1^{p^s},\cdots ,\mathrm{Nm}(\mu)_q^{p^s})\subset  K^{p^s}\]
as ideals in $\bF_p[a_{\spoke},u_{\lambda}^{\pm 1},\mathrm{Nm}(\mu)_1,\cdots,\mathrm{Nm}(\mu)_q]\langle \tilde{\mu}_1,\cdots ,\tilde{\mu}_q,u_{\spoke}\rangle$.
\end{proof} 

\section{The May spectral sequence and the Segal conjecture}\label{mayss}
The goal of this section is to prove \autoref{thm:Segal} using the spectral sequence \eqref{SSa-inverted}. To compute the $\mathrm{E}_2$-term  of the spectral sequence
  \eqref{SSa-inverted}, we use a May spectral sequence: For each integer $n$, filtering by powers of the coaugmentation coideal as in \cite[Example A1.3.10]{Rav86} produces a spectral sequence 
\begin{equation}\label{may-SS} \tag{MSS}
^{(n)}\mathrm{E}_1^{i+k\spoke,s,f}\implies \mathrm{Ext}_{\mathbb{F}_p[\mathrm{Nm}(\mu)]\langle \tilde{\mu}\rangle/(\mathrm{Nm}(\mu)^{p^{n}})}^{s,s+i+k\spoke}(\bF_p,\bF_p[a_{\spoke},u_{\lambda}^{\pm 1}]\langle u_{\spoke}\rangle ) \,. 
\end{equation}
Explicitly, we let $\Gamma= \mathbb{F}_p[\mathrm{Nm}(\mu)]\langle \tilde{\mu}\rangle/(\mathrm{Nm}(\mu)^{p^{n}})$ and consider the short exact sequence 
\[ 
0 \to \mathbb{F}_p \to \Gamma \to \overline{\Gamma} \to 0 
\]
and then consider the increasing filtration given by 
\begin{equation}\label{may-filt}
F_s\Gamma = \mathrm{ker} \left ( \Gamma \overset{\Delta^s}{\longrightarrow} \Gamma^{\otimes s+1} \to  \overline{\Gamma}^{\otimes s+1}  \right )  \,,
\end{equation}
which we call the May filtration. The spectral sequence \eqref{may-SS} is then the spectral sequence associated to the May filtration \eqref{may-filt}. Here we also filter the $\Gamma$-comodule $M=\mathbb{F}_p[a_{\spoke},u_{\lambda}^{\pm 1}]\langle u_{\spoke}\rangle$ by $F_sM=M\otimes_{\Gamma}F_s\Gamma$. 
The only caveat is that this filtration is not compatible with the 
algebra structure on $M$ so the associated May spectral sequence does not have a
Leibniz rule. Nevertheless, it is a module over the May spectral sequence
\[ 
\mathrm{Ext}_{E_0\Gamma}^{*,*,*}(\mathbb{F}_p,\mathbb{F}_p)
\implies 
\mathrm{Ext}_{\Gamma}^{*,*}(\mathbb{F}_p,\mathbb{F}_p) \,,
\]
which has a differential that does have a Leibniz rule. 

Since the spectral sequence takes place in the category of finite type graded
$\mathbb{F}_p$-vector spaces, which is closed under kernels and images, and sequential limits are exact in this category, we can apply the limit to produce a spectral sequence 
\[ 
\mathrm{E}_1^{i+k\spoke,s,f}
\implies  
\lim_n\mathrm{Ext}^{s,s+i+k\spoke}_{\bF_p
[\mathrm{Nm}(\mu)]\langle \tilde{\mu}\rangle/(\mathrm{Nm}(\mu)^{p^{n}})}(\mathbb{F}_p,\mathbb{F}_p[a_{\spoke},u_{\lambda}^{\pm 1}]\langle u_{\spoke} \rangle )  
\]
where
\[
\mathrm{E}_1^{i+k\spoke,s,f}= \lim_n{}^{(n)}\mathrm{E}_1^{i+k\spoke,s,f}\,.
\]

\begin{lem}\label{prop:E1}
We can identify 
\[ 
^{(n)}\mathrm{E}_1^{*,*,*}=\bF_p[a_{\spoke},u_{\lambda}^{\pm 1},x'_{0}, \cdots ,x'_{n-1},z]\langle u_{\spoke} , x_0,\cdots ,x_{n-1} \rangle
\]
where $z$ is represented by $[\tilde{\mu}]$, $x_{j}$ is represented by 
$[\mathrm{Nm}(\mu)^{p^j}]$ and $x'_j$ is represented by 
\[ \sum_{i=1}^{p-1}\frac{1}{p}\binom{p}{i}
[(\mathrm{Nm}(\mu)^{p^j})^{i}|(\mathrm{Nm}(\mu)^{p^j})^{p-i}]  \] 
in the cobar complex. Here 
\begin{align*}
|a_{\spoke}|&=(-\spoke,0, 0) \,,\\ 
|u_{\spoke}|&=(1-\spoke ,0,0) \,,\\ 
|u_{\lambda}|&=(2-2\spoke,0,0) \,,\\ 
|x_j|&=(2p^j-1+ 2p^j(p-1)\spoke,1,1) \,,\\
|x'_j|&=(2p^{j+1}-2+ 2p^{j+1}(p-1)\spoke,2,p) \,, \\
|z| &= (\spoke,1,1) 
\end{align*}
where the first coordinate is the total degree (internal degree minus cohomological degree), the second coordinate is the cohomological degree and the third coordinate is the May filtration. 
\end{lem}

\begin{proof}
Let $\Gamma= \mathbb{F}_p[\mathrm{Nm}(\mu)]\langle \tilde{\mu}\rangle/(\mathrm{Nm}(\mu)^{p^{n}})$ and let the May filtration $F_s\Gamma$ be defined as in \eqref{may-filt}. 
Note that the associated graded of the May filtration is 
\[
E_0(\mathbb{F}_p[\mathrm{Nm}(\mu)]\langle \tilde{\mu}\rangle /(\mathrm{Nm}(\mu)^{p^{n}})= \bigotimes_{i=1}^{n}\mathbb{F}_p[\mathrm{Nm}(\mu)^{p^{i-1}}]/(\mathrm{Nm}(\mu)^{p^i})\langle \tilde{\mu}\rangle  \,.
\]
The computation 
\[
\mathrm{Ext}^{s,s+i+k\spoke}_{E_0(\mathbb{F}_p[\mathrm{Nm}(\mu)]\langle \tilde{\mu}\rangle
  /(\mathrm{Nm}(\mu)^{p^{n}}) }(\mathbb{F}_p,\mathbb{F}_p)=\bF_p[z,x'_{s} :  0 \leq s \leq n-1]\langle x_s : 0 \leq s \leq n-1 \rangle
\]
follows from~\cite[Lemma~3.2.4]{Rav86} for example, where the classes have the cobar representatives given in the statement of the lemma. 
Since we have coefficients in $\bF_p[a_{\spoke}^{\pm 1},u_{\lambda}^{\pm 1}]\langle u_{\spoke} \rangle$
with right coaction given by $\eta_R$ and
\begin{align*} 
(\eta_L-\eta_R)(u_{\spoke})& = 0 \mod F_1\Gamma\,, \\ 
  (\eta_{L}-\eta_R)(u_{\lambda})&=0 \mod F_1\Gamma \,, \text{ and } \\
  (\eta_L-\eta_R)(a_{\spoke})&=0 
\end{align*}
by~\autoref{thm:sTHH}, the result follows. 
\end{proof}
\begin{prop} \label{prop:tate-of-norm}
  There is an isomorphism
\[ 
\pisp(N_e^{C_p}\bF_p)_{a_{\spoke}}^{\wedge}[a_{\spoke}^{-1}]\cong \mathbb{F}_p[a_{\spoke}^{\pm 1}].
\]
\end{prop}
\begin{proof}
We claim that 
\[ 
\mathrm{E}_{\infty}^{*,*,*}\cong\mathbb{F}_p[a_{\spoke}]
\]
in tridegrees contributing to $i+j\spoke$ where $i+j<0$ in the abutment of the spectral sequence \eqref{may-SS}. 
From this we can determine that the $\mathrm{E}_2$-term  of the spectral sequence \eqref{SSa-inverted} is  
\[
\mathrm{E}_2^{i+j\spoke,s}\cong \mathbb{F}_p[a_{\spoke}^{\pm 1}],
\]
since $a_{\spoke}$-towers on any element would eventually land in these bidegrees.
But this is concentrated on the $i=0$ line so it collapses without extensions and 
\[ 
\pisp (N_e^{C_p}\mathbb{F}_p)_{a_{\spoke}}^{\wedge}[a_{\spoke}^{-1}] =\mathbb{F}_p[a_{\spoke}^{\pm 1}]\,. 
\]

We now prove the claim. 
It is clear that in filtration zero we have $\bF_p[a_{\spoke},
u_{\lambda}^{\pm 1}]\langle u_{\Yright} \rangle$. We will see that the $u_{\lambda}$-classes and $u_{\spoke}$-classes support differentials
to kill all element in positive filtrations in degrees $i+j\spoke$ for $i+j< 0$.

Let $I=(i_0,\cdots ,i_{n-1})$ be a tuple of non-negative integers and let $E=(\varepsilon_0,\cdots ,\varepsilon_{n-1})$ be a tuple with $\varepsilon_i\in \{0,1\}$ for each $0\le i\le n-1$. 
Let 
\[
x_{I,E}=(x'_0)^{i_0}\cdot \ldots \cdot (x'_{n-1})^{i_{n-1}}\cdot x_{0}^{\varepsilon_0}\cdot \ldots \cdot  x_{n-1}^{\varepsilon_{n-1}}\,.
\] 
Then the elements
\[ 
X=a_{\Yright}^mu_{\spoke}^{\varepsilon}u_{\lambda}^{\ell}z^{k}x_{I,E}
\]
form an $\mathbb{F}_p$-basis for ${}^{(n)}\mathrm{E}_1$ where $m,k\ge 0$, $\ell\in \mathbb{Z}$, $\varepsilon\in \{0,1\}$, and $I$ and $E$ are tuples as above. 

The May spectral sequence computing
$\mathrm{Ext}_{\Gamma}^{*,*}(\mathbb{F}_p,\mathbb{F}_p)$ is multiplicative and it
collapses so the differentials in our May spectral sequence \eqref{may-SS} are
$z$-linear, $x_{I,E}$-linear and $a_{\Yright}$-linear.

We first compute the $d_1$-differentials. Although \eqref{may-SS} is not
multiplicative, the $d_1$-differential still satisifes the Leibniz rule
as it is the ``linear term'' of the multiplicative formula $\eta_R$.
The formulas \eqref{ulambda} and \eqref{uspoke} imply total
differentials
\begin{align*}  
  d(u_{\lambda}^{p^t})& = \beta^{p^t}  a_{\spoke}^{2p^{t+1}}x_t, \\
    d(u_{\lambda}^{-p^t})& = -\beta^{p^t}  a_{\spoke}^{2p^{t+1}}u_{\lambda}^{-2p^{t}}x_t -
                      \beta^{p^{t+1}}
                      a_{\spoke}^{2p^{t+2}}u_{\lambda}^{-(p+1)p^{t}}x_{t+1}  - \cdots
                      \mod F_2\Gamma ,\\
d(u_{\spoke})& =  \beta'a_{\spoke}^{2}z \,.
\end{align*}
Consequently, if we filter further each vector space of ${}^{(n)} \mathrm{E}_1^{*,*,*}$ by
setting $\mathrm{fil}_x(x_i) = i+1$ and the other elements to be in filtration 0,
the $d_1$-differentials can be computed step by step via $d_{1,x}$ as below: (note that $\beta^p= \beta$)
\begin{align*}
d_{1,0}(u_{\spoke})& =  \beta'a_{\spoke}^{2}z ,\\
d_{1,t+1}(u_{\lambda}^{p^t})& = \beta  a_{\spoke}^{2p^{t+1}}x_t, \\
d_{1,t+1}(u_{\lambda}^{-p^t})& = -\beta  a_{\spoke}^{2p^{t+1}}u_{\lambda}^{-2p^{t}}x_t.
\end{align*}
This gives in total degrees $i+j\spoke$ for $i+j<0$,
\begin{equation*}
^{(n)}\mathrm{E}_2^{*,*,*}=\bF_p[a_{\spoke},u_{\lambda}^{\pm p^n},x'_{0}, \cdots ,x'_{n-1}]\langle u_{\lambda}^{p-1} x_0,\cdots ,u_{\lambda}^{p^{n-1}(p-1)}x_{n-1} \rangle.
\end{equation*}

Moreover, in the cobar complex there are total differentials (note that $\beta^{p-1}=1$)
 \newcommand{\nm}{\mathrm{Nm}(\mu)}
 
\begin{equation*}
  \begin{split}
    &   d\bigg(\sum_{i=1}^{p-1}\frac{1}{p}\binom{p}{i}\beta^{i-1}a_{\spoke}^{2p^{t+1}(i-1)}u_{\lambda}^{p^t(p-i)}[(\nm^{p^{t}})^{i}]\bigg) \\
    & =  a_{\spoke}^{2p^{t+1}(p-1)}  \sum_{i=1}^{p-1}\frac{1}{p}\binom{p}{i}
[(\mathrm{Nm}(\mu)^{p^t})^{i}|(\mathrm{Nm}(\mu)^{p^t})^{p-i}] ,
  \end{split}
\end{equation*}
which gives 
\begin{equation*}
d_{p-1}(u_{\lambda}^{p^t(p-1)} x_{t}) =  a_{\spoke}^{2p^{t+1}(p-1)} x_{t}' \,.
\end{equation*}

Consequently, in total degrees $i+j\spoke$ for $i+j<0$, there is
\begin{equation*}
^{(n)}\mathrm{E}_p^{*,*,*}=\bF_p[a_{\spoke},u_{\lambda}^{\pm p^n}].
\end{equation*}
This proves the claim. 
\end{proof}

\begin{proof}[Proof of \autoref{thm:Segal}] 
By \autoref{prop:tate-of-norm}, we know that 
\[
\pisp (N_e^{C_p}\mathbb{F}_p)_{a_{\spoke}}^{\wedge}[a_{\spoke}^{-1}]=\mathbb{F}_p[a_{\spoke}^{\pm 1}]\,.\] 
We also know that  
\[
\pisp N_e^{C_p}\mathbb{F}_p[a_{\spoke}^{-1}]=
\pi_{*}(N_e^{C_p}\mathbb{F}_p)^{\Phi C_p}[a_{\Yright}^{\pm 1}] = \mathbb{F}_p[a_{\spoke}^{\pm 1}]\,.\] 
So the right arrow in the pullback square 
\[
\begin{tikzcd}
N_e^{C_p}\mathbb{F}_p \ar[r]\ar[d]  & N_e^{C_p}\mathbb{F}_p[a_{\spoke}^{-1}] \ar[d] \\
(N_e^{C_p}\mathbb{F}_p)_{a_{\spoke}}^{\wedge} \ar[r] & (N_e^{C_p}\mathbb{F}_p)_{a_{\spoke}}^{\wedge}[a_{\spoke}^{-1}]
\end{tikzcd}
\]
is an equivalence, and do is the left arrow as desired. 
\end{proof}

\appendix

\section{The \texorpdfstring{$\ROs(C_p)$}{ROCp}-graded homotopy of \texorpdfstring{$\mF_p$}{Fp}}\label{appendixA}

In this appendix, we recall the computation of the
$\ROs(C_p)$-graded homotopy of $\mF_p$.

Let $\lambda_i \in \RO(C_p)$ denote the representation of $C_p$ on $\bR^2$, where a chosen generator of $C_p$ rotates by $e^{2\pi i / p}$.
Working $p$-locally, there is an equivalence $S^{\lambda_i} \simeq S^{\lambda_j}$ \cite[Section
2]{HKSZ22}, so we may take $\lambda = \lambda_1$ and focus on $m + n\lambda \in \RO(C_p)$ for
$m,n \in \bZ$.

D. Wilson~\cite{Wil17} observed that adding the \emph{spoke sphere} $S^{\Yright}$ to
$\RO(C_p)$, $C_p$-equivariant calculations behave nicer and closer to their
$C_2$-equivariant counterparts.
\begin{defn}
In the stable $C_p$-equivariant category, we denote by $S^\Yright$ the cofiber of $(C_p)_+\rightarrow S^0$. We denote by $S^{-\Yright}$ its Spanier--Whitehead dual. 
\end{defn}

There are (non-canonical) equivalences
\begin{equation}
  \label{eq:spoke}
 S^{\lambda} \otimes S^{-\Yright} \simeq S^{\Yright} \text{ and }
S^{\Yright} \otimes S^{-\Yright} \simeq S^0 \oplus \bigoplus_{p-2} (C_p)_+.
\end{equation}

\begin{defn}
Given a $C_p$-spectrum $X$, we denote by $\pisp(X)$ the graded abelian group of
$C_p$-maps $S^{\varoast} \to X$ indexed by gradings $\varoast = V \text{ or } V+\Yright$, for $V \in
\RO(C_p)$. For convenience, we sometimes write $\varoast = V - \Yright$, which
means $\varoast = V - \lambda + \Yright$, or write $\varoast = m + n\Yright$, which means
\[
\varoast = \begin{cases}
             m + \frac{n}{2} \lambda & \text{ for even }n\,;\\
             m + \frac{n-1}{2} \lambda + \Yright & \text{ for odd }n \,.
            \end{cases}
\]
\end{defn}

\begin{rem}
Because the spoke sphere is not a Picard element, spoke graded homotopy does not naturally organize into a ring. It also fails to have suspension isomorphisms for spoke suspensions, and instead there are (non-canonical) isomorphisms 
\begin{equation}
  \label{eq:1}
  \pi^{C_p}_{\star+\Yright}(\Sigma^{\Yright}X) \cong \pi^{C_p}_{\star}(X) \oplus \bigoplus_{p-2} \pi^{C_p}_{\star}((C_p)_+\otimes X).
\end{equation}
\end{rem}

To treat these subtleties, Hahn, Senger, A. Zhang and the second author~\cite{HSZZ} introduced spoke graded
homotopy groups \emph{modulo $a_{\Yright}$-torsion}. We recall this notion below. First we need some definitions. 
\begin{defn}
For $V\in \{\lambda ,\spoke\}$, we define $a_{V}\in \pi_{-V}^{C_p}\mathbb{S}$ to be the
equivalence class of the (desuspension) of the embedding $S^0\to
S^{V}$ of $\{0,\infty\}\in S^{V}$ and also write $a_V\in \pi_{-V}^{C_p}\Ff$
for its Hurewicz image.

We write $u_{\lambda}\in \pi_{2-\lambda}^{C_p}\Ff \cong \mathrm{H}_2^{C_p}(S^\lambda , \mF_p)$ for the class that restricts to a generator of
$\mathrm{H}_{2}(S^2,\mathbb{F}_p) = \mathbb{F}_p$. 
\end{defn}
\begin{rem} \label{rem:a-spoke-square}
 The Spanier--Whitehead dual $Da_{\Yright}: S^{-\Yright} \to S^0$ is $a_{\Yright}$ modulo the
 summands with the free $C_p$'s in \eqref{eq:1}, and the composite $S^0
 \overset{a_{\Yright}}{ \longrightarrow } S^{\Yright} \overset{Da_{\Yright}}{ \longrightarrow } S^{\lambda}$ is $a_{\lambda}$. 
\end{rem}

For a $C_p$-equivariant ring spectrum $X$,  $\varoast$-graded homotopy groups of
$X$ forms a ring modulo $a_{\Yright}$-torsion. In \cite{HSZZ}, this is made precise
with the following \autoref{defn:modtorsion} and \autoref{prop:product}, whose
details we omit. 
\begin{defn}\label{defn:modtorsion}
We define $\pispa(X) $ to be the image of the map 
\[
\pisp(X) \to \pisp(X)[a_{\lambda}^{-1}] \,.
\]
\end{defn}
\noindent Intuitively, $\pispa(X)$ picks out the ``postitive cone'' elements of $\pisp(X)$.
\begin{prop}
\label{prop:product}
For a $C_p$-equivariant ring spectrum $X$,  $\pispa(X)$ forms a
ring. Explicitly, for $a \in \pi^{C_p}_{V + \Yright}(X)$ and $b \in
\pi^{C_p}_{W+\Yright}(X)$, the product  $ab$ is  
\begin{equation*}
ab: S^{V + W + \lambda} \longrightarrow S^{V + \Yright} \otimes S^{W + \Yright} \overset{a \otimes b}{ \longrightarrow } X \otimes X
\longrightarrow X.
\end{equation*}
\end{prop}
Note that $a_{\Yright}^2 = a_{\lambda}$ by \autoref{rem:a-spoke-square}.

\begin{lem}
  \label{lem:u_spoke}
  There exists a class $u_{\spoke}\in \pi_{1-\spoke}^{C_p}\Ff \cong
  \mathrm{H}_{1}^{C_p}(S^{\Yright}, \mF_p)$ that  restricts to $(p-1)$-tuple of generators 
$\mathrm{H}_1(\bigoplus_{p-1} S^1, \mathbb{F}_p)$. Furthermore, $a_{\Yright} u_{\Yright} \neq 0$.
\end{lem}
\begin{proof}
  By calculating the $H\mF_p$-homology long exact sequence
  for the $C_p$-equivariant cofiber sequence $(C_p)_+ \to S^0 \to S^{\Yright}$, we have
  an exact sequence
  $$0 \to \underline{\mathrm{H}}_{1}^{C_p}(S^{\Yright}, \mF_p) \to
  \begin{tikzcd}
    \mathbb{F}_p \ar[d, shift right, "{\triangle}"'] \ar[r,"0"] &   \mathbb{F}_p \ar[d, shift right, "1"'] \\
    \oplus_{p} \mathbb{F}_p \ar[u,shift right, "{\triangledown}"'] \ar[r,"{\triangledown}"] &     \mathbb{F}_p \ar[u,shift right, "0"']
  \end{tikzcd}
   \,.$$ 
  So we conclude 
  \[
  \underline{\mathrm{H}}_{1}^{C_p}(S^{\Yright}, \mF_p) =
  \begin{tikzcd}
    \mathbb{F}_p \ar[d, shift right, "{\triangle}"'] \\
    \oplus_{p-1} \mathbb{F}_p \ar[u,shift right, "0"'] 
  \end{tikzcd}
  \,.\] 
  Let  $u_{\spoke}$ be a generator of $\mathrm{H}_{1}^{C_p}(S^{\Yright}, \mF_p) = \mathbb{F}_p$. For the
  sake of contradiction assume that $a_{\Yright} u_{\Yright} = 0$. Then there exists a
  factorization through the cofiber, as the dotted arrow $f$ in the following diagram:
\begin{equation*}
   \begin{tikzcd}
     S^{-\Yright} \otimes S^{1-\Yright} \ar[rd, "a_{\Yright} u_{\Yright}"]  \ar[d]\\
     S^{1-\Yright} \ar[r,"u_{\Yright}"]
    \ar[d," \mathrm{tr}"] & \Ff\\
     (C_p)_+  \otimes S^{1-\Yright} \ar[ru, dotted, "f"']
  \end{tikzcd}
\end{equation*}
This implies $u_{\Yright} = \mathrm{tr}(f)$.
However, we know $u_{\Yright}$ is not in the image of transfer.
\end{proof}

  \begin{prop} 
    \label{prop:homology-of-spoke}
    The $\RO(C_p)$-graded and $\ROs(C_p)$-graded homotopy
    groups of $\Ff$ are 
    \begin{align*}
      \piro(\Ff) & = \mathbb{F}_p[a_{\lambda}, u_{\lambda}] \langle \kappa_\lambda \rangle \oplus 
    \bZ/p\{\Sigma^{-1}\kappa_{\lambda}^{\varepsilon}u_{\lambda}^{-j}a_{\lambda}^{-k}|\varepsilon \in \{0,1\}, j,k\ge 1
                   \},\\
      \pisp(\Ff) & = \mathbb{F}_p[a_{\Yright},u_{\lambda}]\langle u_{\Yright} \rangle \oplus
    \bZ/p\{\Sigma^{-1}u_{\Yright}^{\varepsilon}u_{\lambda}^{-j}a_{\Yright}^{-k}|\varepsilon \in \{0,1\}, j,k\ge 1  \}.
    \end{align*}
    Here, $\kappa_{\lambda} = a_{\Yright}u_{\Yright}$ and thus $|\kappa_{\lambda}| = 1-\lambda$. The first summand
    is called the positive cone and the second summand is called the negative
    cone.
    
    The $\RO(C_p)$- and $\ROs(C_p)$-graded homotopy
    groups of $\Sigma^{\Yright}\Ff$ are
\begin{align*}
  \piro(\Sigma^{\Yright} \Ff) & \cong  \piro(\Ff) \{\hat{x},\ul{x}\} / \{\kappa_{\lambda} \hat{x} =
                      a_{\lambda} \ul{x},  \kappa_{\lambda} \ul{x} = 0\}, \\
  \pisp(\Sigma^{\Yright} \Ff) & \cong  \pisp(\Ff) \{1^{\Yright}\} \oplus \mathbb{F}_p[u_\lambda^{\pm}]\{\tilde{1}_1, \cdots, \tilde{1}_{p-2}\},
\end{align*}
where $| \hat{x}| = 0$, $| \ul{x}| = 1$, $|1^{\Yright}| = \Yright$ and $|\tilde{1}_i|=\Yright$ for
$1 \leq i \leq p-2$.
\end{prop}
  \begin{rem}
  Writing $\theta =  \Sigma^{-1}\kappa_{\lambda}u_{\lambda}^{-1}a_{\lambda}^{-1} $ with $|\theta| = \lambda-2$, one
  can then write  
\begin{align*}
  \piro(\Ff)_{\negative} & = \{\frac{\theta}{x} | x \in \piro(\Ff)_{\pos}\}, \\
  \pisp(\Ff)_{\negative} & = \{\frac{\theta}{x} | x \in \pisp(\Ff)_{\pos}\}.
\end{align*}
It turns out that this description of the negative cone is a consequence of the Anderson duality and
$\theta$ is the Anderson dual of $1$ \cite{HMZ}.
\end{rem}

\begin{proof}[Proof of \autoref{prop:homology-of-spoke}] 
   (1) We first use the isotropy separation diagram 
   \[
\begin{tikzcd}
(\mathbb{S}^{-V}\otimes \mF_p)_{hC_p} \arrow{r}  \arrow{d} & (\mathbb{S}^{-V}\otimes \mF_p)^{C_p} \arrow{r} \arrow{d}  & (\mathbb{S}^{-V}\otimes \mF_p)^{\Phi C_p}  \arrow{d} \\ 
(\mathbb{S}^{-V}\otimes \mF_p)_{hC_p}\arrow{r} & (\mathbb{S}^{-V}\otimes \mF_p)^{hC_p}  \arrow{r}  & (\mathbb{S}^{-V}\otimes H\mathbb{F}_p)^{tC_p}.
\end{tikzcd}
\]
   to compute $\piro(\Ff)$.
The $\RO(C_p)$-graded homotopy fixed point spectral sequence 
  \begin{align*}
    E_2^{s,t}(V)=H^s(BC_p, \pi_{t+|V|} \Fp) & \implies
    \pi_{V+t-s}\Ff^{hC_p}
  \end{align*}
together with the fact that
  \[ H^*(BC_p, \mathbb{F}_p)= \mathbb{F}_p[\frac{a_{\lambda}}{u_{\lambda}}] \langle \epsilon \rangle \text{ where } |\frac{a_{\lambda}}{u_{\lambda}}| = -2, |\epsilon|=-1\]
gives
\[
\pi_{\star}\Ff^{hC_p} = E_2^{*,*}(\star)=\mathbb{F}_p[a_{\lambda},u_{\lambda}^{\pm 1}]\langle \kappa_{\lambda} \rangle \, \text{ where }\kappa_{\lambda}=\frac{\epsilon a_{\lambda}}{u_{\lambda}}.
\] 
Similarly, the $\RO(C_p)$-graded $C_p$-Tate spectral sequence
\[
\hat{H}^*(BC_p, \pi_{|\star|}\Fp)\implies \pi_{\star}\Ff^{tC_p}
\]
collapses and gives 
\[
\pi_{\star}\Ff^{tC_p} =\mathbb{F}_p[a_{\lambda}^{\pm 1},u_{\lambda}^{\pm 1}]\langle \kappa_{\lambda} \rangle \,.
\] 
The map
\[ \pi_{\star}\Ff^{hC_p}\to \pi_{\star}\Ff^{tC_p}\]
is injective. We conclude that 
\[\pi_{\star} (\Ff)_{hC_p}= \bZ/p\{\Sigma^{-1}\kappa_{\lambda}^{\varepsilon}u_{\lambda}^ja_{\lambda}^{-k} |
  \varepsilon \in \{0,1\}, j\in \mathbb{Z},k\ge 1  \} \,.
\]
Now, the geometric fixed points are
\[ 
\pi_{\star}\Ff^{\Phi C_p}=\pi_*\Ff^{\Phi C_p}[a_{\lambda}^{\pm 1}] = \mathbb{F}_p[a_{\lambda}^{\pm}, u_{\lambda}]\langle \kappa_{\lambda} \rangle
\]
and the map
\[ \pi_{\star}\Ff^{\Phi C_p}\to \pi_{\star}\Ff^{tC_p}\]
is also injective.
Consequently, we determine that 
\begin{align*}
\piro\Ff & = \mathrm{ker}(\pi_{\star}\Ff^{\Phi C_p} \to \Sigma \pi_{\star} (\Ff)_{hC_p}) \oplus
           \mathrm{coker}(\Sigma^{-1}\pi_{\star}\Ff^{\Phi C_p} \to  \pi_{\star}
           (\Ff)_{hC_p})\\
  & =\mathbb{F}_p[a_{\lambda}, u_{\lambda}]\langle \kappa_{\lambda} \rangle \oplus
    \bZ/p\{\Sigma^{-1}\kappa_{\lambda}^{\varepsilon}u_{\lambda}^{-j}a_{\lambda}^{-k}|\varepsilon \in \{0,1\}, j,k\ge 1  \} \,.
\end{align*}
\begin{figure}[!ht]
    \begin{subfigure}{0.35\textwidth}
      \printpage[name = HstarSpoke]
      \caption{$\pi^{C_p}_{\star+\Yright}(\Ff)_{\pos}$}
    \end{subfigure}
    \begin{subfigure}{0.45\textwidth}
      \printpage[name = Hstar]
      \caption{$\piro(\Ff)_{pos}$}
    \end{subfigure}
    \begin{subfigure}{0.45\textwidth}
      \centering
      \printpage[name = HstarFull]
      \caption{$\pisp(\Ff)_{\pos}$}
      \label{fig:pispH}
    \end{subfigure}
    \caption{$\pisp(\Ff)_{\pos}$}
    \label{fig:pispHFull}
  \end{figure}
   \begin{figure}
      \centering
      \printpage[name = HstarFullShift]
      \caption{$\pisp(\Sigma^{-\Yright}\Ff)$}
      \label{fig:pispHShift}
  \end{figure}

(2) Next, we consider the cofiber sequence $S^{-\Yright} \overset{a_{\Yright}}{ \to } S^0 \to (C_p)_+$.
There is an identification
$$\piro((C_p)_+ \otimes \Ff) = \pi_{|\star|}(\Fp) =  \mathbb{F}_p[u_{\lambda}^{\pm}]$$  and the map
$\piro\Ff \to \piro((C_p)_+ \otimes \Ff)$ is the restriction map which sends $1$ to
$1$, $u_{\lambda}$ to $u_{\lambda}$ and $u_{\lambda}^{-1}$ to $0$.
Therefore, we get
\begin{align*}
\piro(\Sigma^{-\Yright} \Ff) & \cong \mathbb{F}_p[a_{\lambda}, u_{\lambda}]\{a_{\lambda}\} \oplus  \mathbb{F}_p[a_{\lambda},
u_{\lambda}]\{\kappa_{\lambda}\} \\
& \oplus \bZ/p\{\Sigma^{-1}\kappa_{\lambda}^{\varepsilon}u_{\lambda}^{-j}a_{\lambda}^{-k}|\varepsilon \in \{0,1\}, j,k\ge 1  \}
    \oplus \bZ/p\{\Sigma^{-1}u_{\lambda}^{-j}|j\ge 1  \},
\end{align*}
where $a_{\lambda} \in \pi_{-\lambda}^{C_p}(\Sigma^{-\Yright} \Ff)$ and $\kappa_{\lambda} \in \pi_{1-\lambda}^{C_p}(\Sigma^{-\Yright} \Ff)$.
This is an expression as a module over $\piro(\Ff)$ (at least for the positive
cone elements).
Moreover, the map $a_{\Yright}: \Sigma^{-\Yright}\Ff \to \Ff$ sends $a_{\lambda}$ to $a_{\lambda}$ and $\kappa_{\lambda}$ to $\kappa_{\lambda}$.
Because $a_{\lambda} = a_{\Yright}^2$ and $\kappa_{\lambda} = a_{\Yright} u_{\Yright}$ up to units
(\autoref{lem:u_spoke}), taking degrees into considerations we may write
$$
\pi^{C_p}_{\star+\Yright}(\Ff)_{\pos}  \cong \piro(\Sigma^{-\Yright} \Ff)_{\pos} \cong  \mathbb{F}_p[a_{\lambda}, u_{\lambda}]\{a_{\Yright}\} \oplus  \mathbb{F}_p[a_{\lambda},
u_{\lambda}]\{u_{\Yright}\}
$$
where $a_{\Yright} \in \pi^{C_p}_{-\lambda+\Yright}(\Ff)$ and $u_{\Yright} \in \pi^{C_p}_{1-\lambda+\Yright}(\Ff)$.

Combining the results, we have $\pisp(\Ff)_{\pos} = \mathbb{F}_p[a_{\Yright},
u_{\lambda}]\{u_{\Yright}\} $. See~\autoref{fig:pispHFull} for an illustration. 

For the module and multiplicative structure on the negative cone, consider
instead the cofiber sequence 
\[
(C_p)_+ \to S^0 \overset{a_{\Yright}}{ \to } S^{\Yright}\,.
\] 
The map $\piro((C_p)_+ \otimes \Ff) \longrightarrow \piro\Ff $
is the transfer map which sends $1$ to
$p=0$, $u_{\lambda}$ to $0$  and $u_{\lambda}^{-j}$ to
$\Sigma^{-1}\kappa_{\lambda}u_{\lambda}^{-j}a_{\lambda}^{-1}$. Arguing similarly as before, we can
deduce that
\[
\pisp(\Ff)_{\negative} = 
\bZ/p\{\Sigma^{-1}u_{\Yright}^{\varepsilon}u_{\lambda}^{-j}a_{\Yright}^{-k}|\varepsilon \in \{0,1\}, j,k\ge 1  \}\,.
\]

(3) The description of $\piro(\Sigma^{\Yright} \Ff) $ follows from the description of
$\piro(\Sigma^{-\Yright} \Ff)$ in part (2), where $\hat{x}$ is the $\lambda$-suspension of $a_{\lambda}$ and
$\ul{x}$ is the $\lambda$-suspension of $\kappa_{\lambda}$. The description of
$\pisp(\Sigma^{\Yright}\Ff)$ can be argued similarly as in part (2), while the classes $\tilde{1}_i$ for
$1 \leq i \leq p-2$ come from the free $C_p$-summands in $S^{\Yright} \otimes S^{-\Yright}$ in the computation
$$\pi_{\Yright}^{C_p}(\Sigma^{\Yright}\Ff) \cong [S^{\Yright}, \Sigma^{\Yright}\Ff]^{C_p} \cong [S^0 \oplus \bigoplus_{p-2}(C_p)_+, \Ff]^{C_p}.$$
\end{proof}
The following corollary immediately follows. 
\begin{cor}
We can identify 
\begin{align*}
\pispa H\mF_p& \cong \mathbb{F}_p[a_{\spoke},u_{\lambda}]\langle u_{\spoke} \rangle  \,,\\ 
\pisp H\mF_p[a_{\spoke}^{-1}]& \cong \mathbb{F}_p[a_{\spoke}^{\pm 1},u_{\lambda}]\langle u_{\spoke} \rangle \text{ and }\\
\pisp (H\mF_p)_{a_{\spoke}}^{\wedge}[a_{\spoke}^{-1}]& \cong \mathbb{F}_p[a_{\spoke}^{\pm 1},u_{\lambda}^{\pm 1}]\langle u_{\spoke}\rangle 
\end{align*}
as $\ROs(C_p)$-graded rings.
\end{cor}

\bibliographystyle{amsalpha}
\bibliography{ref}

@incollection {AlmkvistFossum,
    AUTHOR = {Almkvist, Gert and Fossum, Robert},
     TITLE = {Decomposition of exterior and symmetric powers of
              indecomposable {${\bf Z}/p{\bf Z}$}-modules in characteristic
              {$p$}\ and relations to invariants},
 BOOKTITLE = {S\'eminaire d'{A}lg\`ebre {P}aul {D}ubreil, 30\`eme ann\'ee
              ({P}aris, 1976--1977)},
    SERIES = {Lecture Notes in Math.},
    VOLUME = {Vol. 641},
     PAGES = {1--111},
 PUBLISHER = {Springer, Berlin-New York},
      YEAR = {1978},
      ISBN = {3-540-08665-X},
   MRCLASS = {14L30 (13H10 16A26 20C05)},
  MRNUMBER = {499459},
MRREVIEWER = {H.\ R\"ohrl},
}

@Misc{HMZ,
author={Eric Hogle and Clover May and Foling Zou},
title={Equivariant {E}ilenberg--{M}acLane spectra are self-injective},
 HowPublished = {in preparation},
}

@article {LRZLoday,
    AUTHOR = {Lindenstrauss, Ayelet and Richter, Birgit and Zou, Foling},
     TITLE = {Loday constructions of {T}ambara functors},
   JOURNAL = {J. Algebra},
  FJOURNAL = {Journal of Algebra},
    VOLUME = {683},
      YEAR = {2025},
     PAGES = {278--306},
      ISSN = {0021-8693,1090-266X},
   MRCLASS = {55P91 (13D03)},
  MRNUMBER = {4929908},
       DOI = {10.1016/j.jalgebra.2025.06.016},
       URL = {https://doi.org/10.1016/j.jalgebra.2025.06.016},
}

@article{Goo85,
 author = {Goodwillie, Thomas G.},
 title = {Cyclic homology, derivations, and the free loopspace},
 fjournal = {Topology},
 journal = {Topology},
 issn = {0040-9383},
 volume = {24},
 pages = {187--215},
 year = {1985},
 language = {English},
 doi = {10.1016/0040-9383(85)90055-2},
 zbMATH = {3908625},
 Zbl = {0569.16021}
}

@ARTICLE{AKMP24,
       author = {{Angelini-Knoll}, Gabriel and {Merling}, Mona and {P{\'e}roux}, Maximilien},
        title = "{Topological $\Delta G$ homology of rings with twisted $G$-action}",
      journal = {arXiv e-prints},
         year = 2024,
        month = sep,
          eid = {arXiv:2409.18187},
        pages = {arXiv:2409.18187},
          doi = {10.48550/arXiv.2409.18187},
archivePrefix = {arXiv},
       eprint = {2409.18187},
 primaryClass = {math.AT},
       adsurl = {https://ui.adsabs.harvard.edu/abs/2024arXiv240918187A}
}

@article{DMPR24,
 author = {Dotto, Emanuele and Moi, Kristian and Patchkoria, Irakli},
 title = {On the geometric fixed points of real topological cyclic homology},
 fjournal = {Journal of the London Mathematical Society. Second Series},
 journal = {J. Lond. Math. Soc., II. Ser.},
 issn = {0024-6107},
 volume = {109},
 number = {2},
 pages = {68},
 note = {Id/No e12862},
 year = {2024},
 language = {English},
 doi = {10.1112/jlms.12862},
 zbMATH = {7811233},
 Zbl = {1534.19001}
}

@misc{Hog16,
 author = {H{\o}genhaven, Amalie},
 title = {Real topological cyclic homology of spherical group rings},
 year = {2016},
 howpublished = {Preprint, {arXiv}:1611.01204 [math.{AT}] (2016)},
 url = {https://arxiv.org/abs/1611.01204},
 arXiv = {arXiv:1611.01204}
}

@Misc{HSZZ,
author={Hahn, Jeremy and Senger, Andrew and Zhang, Adela (YiYu) and Zou, Foling},
title={An equivariant Adams spectral sequence for {$\mathrm{tmf}(2)$}},
 HowPublished = {in preparation},
}

@article{BHM89,
 author = {B{\"o}kstedt, M. and Hsiang, W.-C. and Madsen, I.},
 title = {The cyclotomic trace and the {K}-theoretic analogue of {Novikov}'s conjecture},
 fjournal = {Proceedings of the National Academy of Sciences of the United States of America},
 journal = {Proc. Natl. Acad. Sci. USA},
 issn = {0027-8424},
 volume = {86},
 number = {22},
 pages = {8607--8609},
 year = {1989},
 language = {English},
 doi = {10.1073/pnas.86.22.8607},
 url = {europepmc.org/articles/pmc298335},
 zbMATH = {4121024},
 Zbl = {0684.55013}
}

@misc{HHRodd,
      title={On the 3-primary {A}rf-{K}ervaire invariant problem}, 
      author={Hill, M. A. and Hopkins, M. J. and Ravenel, D. C.},
      year={2011},
      note={https://people.math.rochester.edu/faculty/doug/mypapers/odd.pdf},
}

@book{BH21,
 author = {Bachmann, Tom and Hoyois, Marc},
 title = {Norms in motivic homotopy theory},
 fseries = {Ast{\'e}risque},
 series = {Ast{\'e}risque},
 issn = {0303-1179},
 volume = {425},
 isbn = {978-2-85629-939-5},
 year = {2021},
 publisher = {Paris: Soci{\'e}t{\'e} Math{\'e}matique de France (SMF)},
 language = {English},
 doi = {10.24033/ast.1147},
 zbMATH = {7403459},
 Zbl = {1522.14028}
}

@Misc{AKBJK,
author={Gabriel Angelini-Knoll and Mark Behrens and Maxwell Johnson and Hana Jia Kong},
title={A {$C_3$}-equivariant Snaith construction},
 HowPublished = {in preparation},
}

@Misc{AKBJK2,
author={Gabriel Angelini-Knoll and Mark Behrens and Maxwell Johnson and Hana Jia Kong},
title={A delocalized construction of $\mathrm{BP}_{\mu_p}$},
 HowPublished = {in preparation},
}

@misc{AKKQ25,
 author = {Angelini-Knoll, Gabriel and Kong, Hana Jia and Quigley, J. D.},
 title = {Real syntomic cohomology},
 year = {2025},
 howpublished = {Preprint, {arXiv}:2505.24734 [math.{AT}] (2025)},
 url = {https://arxiv.org/abs/2505.24734},
 arXiv = {arXiv:2505.24734}
}

@article{BH15,
 author = {Blumberg, Andrew J. and Hill, Michael A.},
 title = {Operadic multiplications in equivariant spectra, norms, and transfers},
 fjournal = {Advances in Mathematics},
 journal = {Adv. Math.},
 issn = {0001-8708},
 volume = {285},
 pages = {658--708},
 year = {2015},
 language = {English},
 doi = {10.1016/j.aim.2015.07.013},
 zbMATH = {6499680},
 Zbl = {1329.55012}
}

@Misc{Bok85,
 Author = {B{\"o}kstedt, Marcel},
 Title = {The topological {H}ochschild homology of $\mathbb{Z}$ and $\mathbb{Z}/p$},
 Year = {1985},
 HowPublished = {Preprint},
 URL = {https://people.math.rochester.edu/faculty/doug/otherpapers/bokstedt2.pdf},
}

@Article{Bou79,
 Author = {Bousfield, A. K.},
 Title = {The localization of spectra with respect to homology},
 FJournal = {Topology},
 Journal = {Topology},
 ISSN = {0040-9383},
 Volume = {18},
 Pages = {257--281},
 Year = {1979},
 Language = {English},
 DOI = {10.1016/0040-9383(79)90018-1},
 zbMATH = {3649651},
 Zbl = {0417.55007}
}

@Article{DMPR21,
 Author = {Dotto, Emanuele and Moi, Kristian and Patchkoria, Irakli and Reeh, Sune Precht},
 Title = {Real topological {Hochschild} homology},
 FJournal = {Journal of the European Mathematical Society (JEMS)},
 Journal = {J. Eur. Math. Soc. (JEMS)},
 ISSN = {1435-9855},
 Volume = {23},
 Number = {1},
 Pages = {63--152},
 Year = {2021},
 Language = {English},
 DOI = {10.4171/JEMS/1007},
 zbMATH = {7328108},
 Zbl = {1473.16005}
}

@article{GIKR22,
 author = {Gheorghe, Bogdan and Isaksen, Daniel C. and Krause, Achim and Ricka, Nicolas},
 title = {{{\(\mathbb{C}\)}}-motivic modular forms},
 fjournal = {Journal of the European Mathematical Society (JEMS)},
 journal = {J. Eur. Math. Soc. (JEMS)},
 issn = {1435-9855},
 volume = {24},
 number = {10},
 pages = {3597--3628},
 year = {2022},
 language = {English},
 doi = {10.4171/JEMS/1171},
 zbMATH = {7547822},
 Zbl = {1498.14050}
}

@misc{Gun81,
  title        = {Cohomotopy of some classifying spaces},
  author       = {J. H. Gunawardena},
  year         = 1980,
  month        = {},
  address      = {Cambridge},
  note        = {J. T. Knight Prize Essay},
  school       = {},
  type         = {}
}

@Article{HSW23,
 Author = {Hahn, Jeremy and Senger, Andrew and Wilson, Dylan},
 Title = {Odd primary analogs of real orientations},
 FJournal = {Geometry \& Topology},
 Journal = {Geom. Topol.},
 ISSN = {1465-3060},
 Volume = {27},
 Number = {1},
 Pages = {87--129},
 Year = {2023},
 Language = {English},
 DOI = {10.2140/gt.2023.27.87},
 zbMATH = {7688325},
 Zbl = {1523.55011}
}

@Article{HW21,
 Author = {Hahn, Jeremy and Wilson, Dylan},
 Title = {Real topological {Hochschild} homology and the {Segal} conjecture},
 FJournal = {Advances in Mathematics},
 Journal = {Adv. Math.},
 ISSN = {0001-8708},
 Volume = {387},
 Pages = {17},
 Note = {Id/No 107839},
 Year = {2021},
 Language = {English},
 DOI = {10.1016/j.aim.2021.107839},
 zbMATH = {7369662},
 Zbl = {1472.55011}
}

@Article{HW22,
 Author = {Hahn, Jeremy and Wilson, Dylan},
 Title = {Redshift and multiplication for truncated {Brown}-{Peterson} spectra},
 FJournal = {Annals of Mathematics. Second Series},
 Journal = {Ann. Math. (2)},
 ISSN = {0003-486X},
 Volume = {196},
 Number = {3},
 Pages = {1277--1351},
 Year = {2022},
 Language = {English},
 DOI = {10.4007/annals.2022.196.3.6},
 zbMATH = {7611907},
 Zbl = {1541.55010}
}

@InCollection{HN20,
 Author = {Hesselholt, Lars and Nikolaus, Thomas},
 Title = {Topological cyclic homology},
 BookTitle = {Handbook of homotopy theory},
 ISBN = {978-0-8153-6970-7; 978-1-032-91738-2; 978-1-351-25162-4},
 Pages = {619--656},
 Year = {2020},
 Publisher = {Boca Raton, FL: CRC Press},
 Language = {English},
 DOI = {10.1201/9781351251624-15},
 URL = {hdl.handle.net/21.11116/0000-0008-00AA-8},
 zbMATH = {7303328},
 Zbl = {1473.14038}
}

@Article{HHR16,
 Author = {Hill, M. A. and Hopkins, M. J. and Ravenel, D. C.},
 Title = {On the nonexistence of elements of {Kervaire} invariant one},
 FJournal = {Annals of Mathematics. Second Series},
 Journal = {Ann. Math. (2)},
 ISSN = {0003-486X},
 Volume = {184},
 Number = {1},
 Pages = {1--262},
 Year = {2016},
 Language = {English},
 DOI = {10.4007/annals.2016.184.1.1},
 zbMATH = {6605831},
 Zbl = {1366.55007}
}

@Misc{HKSZ22,
 Author = {Hu, Po and Kriz, Igor and Somberg, Petr and Zou, Foling},
 Title = {The {$\mathbb{Z}/p$}-equivariant dual {Steenrod} algebra for an odd prime $p$},
 Year = {2022},
 HowPublished = {Preprint, {arXiv}:2205.13427 [math.{AT}] (2022)},
 URL = {https://arxiv.org/abs/2205.13427},
 arXiv = {arXiv:2205.13427}
}

@article{MSV97,
 author = {McClure, J. and Schw\"anzl, R. and Vogt, R.},
 title = {$\mathrm{THH}(\mathrm{R}) \cong  \mathrm{R} \otimes \mathrm{S}^1$ for $\mathrm{E}_{\infty}$ ring spectra},
 fjournal = {Journal of Pure and Applied Algebra},
 journal = {J. Pure Appl. Algebra},
 issn = {0022-4049},
 volume = {121},
 number = {2},
 pages = {137--159},
 year = {1997},
 language = {English},
 doi = {10.1016/S0022-4049(97)00118-7},
 zbMATH = {1080528},
 Zbl = {0885.55004}
}

@ARTICLE{NS22,
       author = {{Nardin}, Denis and {Shah}, Jay},
        title = "{Parametrized and equivariant higher algebra}",
      journal = {arXiv e-prints},
         year = 2022,
        month = feb,
          eid = {arXiv:2203.00072},
        pages = {arXiv:2203.00072},
          doi = {10.48550/arXiv.2203.00072},
archivePrefix = {arXiv},
       eprint = {2203.00072},
 primaryClass = {math.AT},
       adsurl = {https://ui.adsabs.harvard.edu/abs/2022arXiv220300072N}
}

@article {NS18,
    AUTHOR = {Nikolaus, Thomas and Scholze, Peter},
     TITLE = {On topological cyclic homology},
   JOURNAL = {Acta Math.},
  FJOURNAL = {Acta Mathematica},
    VOLUME = {221},
      YEAR = {2018},
    NUMBER = {2},
     PAGES = {203--409},
      ISSN = {0001-5962,1871-2509},
   MRCLASS = {55U35 (16E40 18E30 19D99)},
  MRNUMBER = {3904731},
MRREVIEWER = {Geoffrey\ M. L. Powell},
       DOI = {10.4310/ACTA.2018.v221.n2.a1},
       URL = {https://doi.org/10.4310/ACTA.2018.v221.n2.a1},
}

@ARTICLE{LLP25,
       author = {{Lenz}, Tobias and {Linskens}, Sil and {P{\"u}tzst{\"u}ck}, Phil},
        title = "{Norms in equivariant homotopy theory}",
      journal = {arXiv e-prints},
         year = 2025,
        month = mar,
          eid = {arXiv:2503.02839},
        pages = {arXiv:2503.02839},
          doi = {10.48550/arXiv.2503.02839},
archivePrefix = {arXiv},
       eprint = {2503.02839},
 primaryClass = {math.AT},
       adsurl = {https://ui.adsabs.harvard.edu/abs/2025arXiv250302839L}
}

@article{LNR12,
 author = {Lun{\o}e-Nielsen, Sverre and Rognes, John},
 title = {The topological {Singer} construction},
 fjournal = {Documenta Mathematica},
 journal = {Doc. Math.},
 issn = {1431-0635},
 volume = {17},
 pages = {861--909},
 year = {2012},
 language = {English},
 doi = {10.4171/dm/384},
 zbMATH = {6155120},
 Zbl = {1275.55005}
}

@Article{Man24,
 Author = {Mantovani, Lorenzo},
 Title = {Localizations and completions of stable {{\(\infty\)}}-categories},
 FJournal = {Rendiconti del Seminario Matematico della Universit{\`a} di Padova},
 Journal = {Rend. Semin. Mat. Univ. Padova},
 ISSN = {0041-8994},
 Volume = {151},
 Pages = {1--62},
 Year = {2024},
 Language = {English},
 DOI = {10.4171/RSMUP/122},
 zbMATH = {7828274},
 Zbl = {1535.18040}
}

@book{Rav86,
 author = {Ravenel, Douglas C.},
 title = {Complex cobordism and stable homotopy groups of spheres},
 fseries = {Pure and Applied Mathematics (Academic Press)},
 series = {Pure Appl. Math., Academic Press},
 issn = {0079-8169},
 volume = {121},
 year = {1986},
 publisher = {Academic Press, New York, NY},
 language = {English},
 zbMATH = {3984103},
 Zbl = {0608.55001}
}

@Article{SW22,
 Author = {Sankar, Krishanu and Wilson, Dylan},
 Title = {On the {{\(C_p\)}}-equivariant dual {Steenrod} algebra},
 FJournal = {Proceedings of the American Mathematical Society},
 Journal = {Proc. Am. Math. Soc.},
 ISSN = {0002-9939},
 Volume = {150},
 Number = {8},
 Pages = {3635--3647},
 Year = {2022},
 Language = {English},
 DOI = {10.1090/proc/15846},
 zbMATH = {7574458},
 Zbl = {1530.55008}
}

@article{Isa01,
 author = {Isaksen, Daniel C.},
 title = {A model structure on the category of pro-simplicial sets},
 fjournal = {Transactions of the American Mathematical Society},
 journal = {Trans. Am. Math. Soc.},
 issn = {0002-9947},
 volume = {353},
 number = {7},
 pages = {2805--2841},
 year = {2001},
 language = {English},
 doi = {10.1090/S0002-9947-01-02722-2},
 zbMATH = {1597939},
 Zbl = {0978.55014}
}

@misc{QS21,
 author = {Quigley, J. D. and Shah, Jay},
 title = {On the equivalence of two theories of real cyclotomic spectra},
 year = {2021},
 howpublished = {Preprint, {arXiv}:2112.07462 [math.{AT}] (2021)},
 url = {https://arxiv.org/abs/2112.07462},
 arXiv = {arXiv:2112.07462}
}

@ARTICLE{Ste25,
       author = {{Stewart}, Natalie},
        title = "{On tensor products with equivariant commutative operads}",
      journal = {arXiv e-prints},
         year = 2025,
        month = apr,
          eid = {arXiv:2504.02143},
        pages = {arXiv:2504.02143},
          doi = {10.48550/arXiv.2504.02143},
archivePrefix = {arXiv},
       eprint = {2504.02143},
 primaryClass = {math.AT},
       adsurl = {https://ui.adsabs.harvard.edu/abs/2025arXiv250402143S}
}

@phdthesis{Wil17,
  title        = {Equivariant, parameterized, and chromatic homotopy theory},
  author       = {Dylan Wilson},
  year         = 2017,
  month        = {June},
  address      = {},
  note         = {},
  school       = {Northwestern University},
  type         = {PhD thesis}
}

@Misc{HKSZ24,
 Author = {Hu, Po and Kriz, Igor and Somberg, Petr and Zou, Foling},
 Title = {The {$\mathbb{Z}/p$}-equivariant spectrum $BP\mathbb{R}$ for an odd prime $p$},
 Year = {2024},
 HowPublished = {Preprint, {arXiv}:2407.16599 [math.{AT}] (2024)},
 URL = {https://arxiv.org/abs/2407.16599},
 arXiv = {arXiv:2407.16599}
}

@article{Sha23,
 author = {Shah, Jay},
 title = {Parametrized higher category theory},
 fjournal = {Algebraic \& Geometric Topology},
 journal = {Algebr. Geom. Topol.},
 issn = {1472-2747},
 volume = {23},
 number = {2},
 pages = {509--644},
 year = {2023},
 language = {English},
 doi = {10.2140/agt.2023.23.509},
 zbMATH = {7706485},
 Zbl = {1541.55020}
}

\end{document}